\documentclass{m2an}
\usepackage{amsmath,amsfonts}
\usepackage{amssymb,mathrsfs,graphicx,bbm}
\usepackage{verbatim,latexsym,caption,tikz,multirow, url, color}
\usepackage{algorithmic}

\newcommand\tcb{\textcolor{blue}}

\newtheorem{assumption}{Assumption}[section]

\numberwithin{equation}{section}

\def\R{\mathbb{R}}
\def\P{\mathbb{P}}
\def\E{\mathbb{E}}

\def\cO{\mathcal{O}}
\def\cL{\mathcal{L}}
\def\cP{\mathcal{P}}
\def\cW{\mathcal{W}}
\def\cC{\mathcal{C}}
\def\cB{\mathcal{B}}
\def\cF{\mathcal{F}}

\def\bfX{\mathbf{X}}

\DeclareMathOperator{\var}{var}
\DeclareMathOperator{\tr}{tr}

\begin{document}
%%-----------------------------
%%      the top matter
%%-----------------------------
\title{Convergence analysis of an explicit method and its random batch approximation for the McKean-Vlasov equations with non-globally Lipschitz conditions}
\author{Qian Guo$^{1}$, Jie He$^{*,}$}\address{Department of Mathematics, Shanghai Normal University, Shanghai 200234, China.}%\thanks{Corresponding author.
%\\ E-mail addresses:{qguo@shnu.edu.cn(Q. Guo), 1000497941@smail.shnu.edu.cn(J. He), leili2010@sjtu.edu.cn (L. Li).}}}
\author{Lei Li}\address{School of Mathematical Sciences, Institute of Natural Sciences, MOE-LSC, DCI joint team, Shanghai Jiao Tong University, 200240, Shanghai, PR China.\\\thanks{Corresponding author.
\\ E-mail addresses:{qguo@shnu.edu.cn(Q. Guo), 1000497941@smail.shnu.edu.cn(J. He), leili2010@sjtu.edu.cn (L. Li).}}}
%\author{...}\address{...}
%
%\date{...}
%
\begin{abstract} In this paper, we present a numerical approach to solve the McKean-Vlasov equations, which are distribution-dependent stochastic differential equations, under some non-globally Lipschitz conditions for both the drift and diffusion coefficients. We establish a propagation of chaos result, based on which the McKean-Vlasov equation is approximated by an interacting particle system. A truncated Euler scheme is then proposed for the interacting particle system allowing for a Khasminskii-type condition on the coefficients. To reduce the computational cost, the random batch approximation proposed in [Jin et al., J. Comput. Phys., 400(1), 2020] is extended to the interacting particle system where the interaction could take place in the diffusion term. An almost half order of convergence is proved in $L^p$ sense.  Numerical tests are performed to verify the theoretical results. \end{abstract}
%
%\begin{resume} In this paper, we present a numerical approach to solve the McKean-Vlasov equations, which are distribution-dependent stochastic differential equations, under some non-globally Lipschitz conditions for both the drift and diffusion coefficients. We establish a propagation of chaos result, based on which the McKean-Vlasov equation is approximated by an interacting particle system. A truncated Euler scheme is then proposed for the interacting particle system allowing for a Khasminskii-type condition on the coefficients. To reduce the computational cost, the random batch approximation proposed in [Jin et al., J. Comput. Phys., 400(1), 2020] is extended to the interacting particle system where the interaction could take place in the diffusion term. An almost half order of convergence is proved in $L^p$ sense.  Numerical tests are performed to verify the theoretical results. \end{resume}
%
%\subjclass{???, ???}
%
\maketitle
%%-----------------------------
%%      your text
%%-----------------------------
\section*{Introduction}
The McKean-Vlasov equations are some distribution dependent stochastic differential equations (SDEs) that are first introduced by H.\ McKean \cite{mckean1966} to describe the Vlasov system in dynamics theory, and mean field stochastic differential equations.
They have been widely studied due to their important applications in nonlinear equations and finance \cite{pham2016,zhang2021,borkar2010}.

In this paper we consider the following McKean-Vlasov SDE:
\begin{equation}\label{McKeanLimit}
X_{t}=X_{0}+\int_{0}^{t}a\left(X_{s}, \mathcal{L}\left(X_{s}\right)\right) d s+\int_{0}^{t} b\left(X_{s}, \mathcal{L}\left(X_{s}\right)\right) d W_{s},%\quad \forall t \in [0, T],
 \end{equation}
 where $\cL(X_s)$ indicates the law of $X_s$, i.e., $\mu_s:=\cL(X_s)$ satisfies
 \[
 \mu_s(E)=\P(X_s\in E)
 \]
with any Borel $E\subset\R^d$.
Here, $X_{0} \in \R^d$ is the initial value of the SDE sampled from an initial distribution $\mu_0$, $a: \R^d\times \cP_{2}(\R^d) \to \R^d$, $b: \R^d\times \cP_{2}(\R^d) \to \R^{d\times m'}$ are some given fields and $W_{\cdot}$ is an $m'$-dimensional standard Wiener processes so $W_s$ is the process at time $s$. Here, $\cP_{2}$ is defined as follows.  We denote  the set of all probability measures on $\R^d$ by $\cP(\R^d)$ while the notation $\cP_p$ indexed by some $p\ge 1$ means the set of probability measures having finite $p$-moment, i.e.,
\begin{equation}\label{eq:P2}
\mathcal{P}_{p}\left(\mathbb{R}^{d}\right):=\left\{\mu \in \mathcal{P}\left(\mathbb{R}^{d}\right): \int_{\mathbb{R}^{d}}|x|^p \mu(d x)<\infty\right\}.
\end{equation}
Clearly, $\cP_p\subset \cP_2$ for $p\ge 2$. We consider the dynamics for all $t\in (0,\infty)$. Later for numerical approximation, we will fix one, but arbitrary, time $T$ to study the convergence. 
Different from the classical SDEs \cite{mao2007}, the coefficients of McKean–Vlasov SDEs depend on the law of the current variable $X_s$ so the dynamics of the law $\mathcal{L}\left(X\right)$ is nonlinear \cite{Kolokoltsov2010}.

The existence and uniqueness for strong  solutions of \eqref{McKeanLimit} have been established in \cite{Sznitman1991} under some linear growth and Lyapunov-type conditions. Later the existence and uniqueness for weak and strong solutions of the McKean–Vlasov SDEs have been studied by some authors under various conditions and we refer the reader to \cite{Bauer2018,Mishura2020,haji2021simple} for more details.
To our knowledge, there is little literature that reports the existence and uniqueness results on the solution of a McKean–Vlasov SDE with  super-linear coefficients. In \cite{dos2019}, the authors  showed that the McKean–Vlasov SDE has a unique solution when the drift coefficient satisfies a one-sided Lipschitz condition, but the diffusion term is still linearly growing. It is worthwhile to mention that a recent result that proves that the McKean-Vlasov SDE admits a unique solution when both drift and diffusion coefficients are super-linear w.r.t. the state \cite{stockinger2021}.

As is well known, the first step to simulate the  McKean--Vlasov SDEs is usually to approximate the true measure $ \mathcal{L}({X_t})$ by the empirical measure. That is, a particle system is adopted to simulate the McKean--Vlasov SDE following the theory of propagation of chaos \cite{Sznitman1991}.  The state of the particle $i\in {1, \ldots, N}$ in the symmetric system of SDEs coupled in a mean field scaling is then given by
\begin{equation}\label{IN}
X_{t}^{i}=X_{0}^{i}+\int_{0}^{t} a\left( X_{s}^{i}, \mu_{s}^{X}\right) d s+\int_{0}^{t} b\left(X_{s}^{i}, \mu_{s}^{X}\right) d W_{s}^{i},
  \end{equation}
where
\[
\mu_{t}^{X}(\cdot):=\frac{1}{N} \sum_{i=1}^{N} \delta_{X_{t}^{i}}(\cdot),
\]
 $\delta_{x}$ denotes the Dirac measure at point $x$,  $X^{i}_{0} \in \R^d$ are the initial values of the SDEs and the processes $W_s^i$ ($1\le i\le N$) are $N$ independent copies of the $m'$-dimensional standard Wiener processes.
We always assume that $X_0^i$ are i.i.d. sampled from some initial law $\mu_0$.

Assuming the coefficients of (\ref{McKeanLimit}) satisfy global Lipschitz condition w.r.t. the state and measure, the convergence of an Euler scheme applied to the particle system has been shown in \cite{Bossy1997}.  In \cite{haji2018multilevel}, multilevel and multi-index Monte Carlo methods have been proposed for the Euler schemes applied to the McKean-Vlasov equation with global Lipschitz conditions. Bao and Huang \cite{bao2021} studied propagation of chaos and convergence rate of the tamed Euler–Maruyama scheme for McKean–Vlasov SDEs by using associated weakly interacting particle systems,  where the drift or diffusion term is  Hölder continuous. Then the convergence of Milstein schemes, taming the drift by a one-sided Lipschitz
condition, for a time-discretized interacting particle system was discussed in \cite{bao2021first}. Li et al. \cite{LiMaoSongWuYin2022} proved strong convergence of the Euler–Maruyama schemes for approximating McKean–Vlasov SDEs under local Lipschitz conditions w.r.t. the
state variable.
Under a Khasminskii-type monotonicity condition, instead of imposing a one-sided and global Lipschitz condition on the drift and diffusion coefficient, respectively, two explicit schemes were proposed in \cite{stockinger2021} for an interacting particle system and the propagation of chaos property of associated McKean–Vlasov equation is also studied therein.  In this paper, we extend the theory of propagation of chaos for a highly nonlinear McKean–Vlasov equation \cite{stockinger2021} to the case of $p\geq 2$.
 Hoeksema et al. \cite{hoeksema} proved a large deviation principle for the empirical measure of a general system of mean-field interacting diffusions with a singular drift and show convergence to the associated McKean–Vlasov equation. Rached et al. introduced in \cite{rached2022single} a decoupling approach for McKean-Vlasov SDEs, which approximates theMcKean-Vlasov law in the coefficients.

Clearly, the computational cost for the numerical schemes is $\cO(N^2)$ per time step if naively implemented, which is expensive in applications.
 To save the computational cost, a projection-based particle method was proposed in \cite{belomestny2018} for solving McKean–Vlasov stochastic differential equations. Recently, a random batch method has been introduced in \cite{jin2020} to simulate the interacting particle systems efficiently, which can reduce the computational cost from $\cO(N^2)$ per time step to $\cO(N)$ with additive noise.

The main contribution of this paper is twofold. Firstly, we provide a strong convergence analysis of a truncated Euler scheme, which is different from the tamed scheme in \cite{stockinger2021},  to numerically solve the McKean-Vlasov SDEs under a Khasminskii-type monotonicity condition w.r.t. the state. Secondly, we extend the aforementioned random batch approximation to improve the efficiency of the truncated method, and establish the error estimate for the approximation. Different from the results in \cite{jin2020}, the diffusion coefficient in the interacting particle systems considered in this paper could contain weakly interacting terms so that the random batch approximation also takes place in the diffusion.

The remainder of this article is organized as follows: In Section \ref{part2}, some assumptions are introduced and the propagation of chaos property is obtained. In Section \ref{part3}, we employ a truncated Euler method to solve the interacting particle systems and prove that the convergence order is one-half. Then in Section \ref{part4}, we present the convergence and efficiency of the numerical method by integrating the random batch method into the truncated Euler scheme. In the last section (Section \ref{part5}), a numerical experiment is given to verify the theoretical results.

\section{Mathematical preliminaries}\label{part2}

We first introduce some notations and basic definitions for later sections. The Euclidean norm of a $d$-dimensional vector is denoted by $|\cdot|$ and the Hilbert-Schmidt norm of a $d\times m$-matrix is denoted by $||\cdot||$, where we recall the Hilbert-Schmidt norm is given by $\|A\|^2=\tr(A^TA)$ for a real matrix $A$. If $G $ is a set, its indicator function $\mathcal{I}_ {G} $ is given by  $$\mathcal{I}_ {G} =\begin{cases} 1,  \quad x \in G, \\0, \quad\text{otherwise}.\end{cases}$$ %then $\mathcal{I}_ {G} =1 $, otherwise 0.} %\mathcal{I}_{G}
Moreover, $a\wedge b :=\min(a, b)$ and $a\vee b :=\max(a, b)$. The notation $u\otimes v$ for $u\in \R^d$ and $v\in \R^d$ means the tensor product of $u$ and $v$, namely, $(u\otimes v)_{ij}=u_i v_j$. Throughout this article, $C > 0$ is a generic constant that might change its value from line to line. We use the notation $A \lesssim B $ to mean $A\le C B$ for some generic constant $C$ that is independent of the parameters been focused on.

 Recall \eqref{eq:P2} for $\cP_p$.  For $p\ge 1$ and any $\mu, \nu\in\mathcal{P}_{p}(\mathbb{R}^{d})$, the $p$-Wasserstein distance is defined by
 \begin{equation*}
 \mathcal{W}_{p}(\mu, \nu):=\left(\inf _{\pi \in \Pi(\mu, \nu)} \int_{\mathbb{R}^{d} \times \mathbb{R}^{d}}|x-y|^{p} \pi(d x, d y)\right)^{1/p},
 \end{equation*}
where $\Pi(\mu, \nu)$ is the set of couplings of $\mu$ and $\nu$ (i.e., joint distributions with marginals to be $\mu$ and $\nu$ respectively). Clearly, $\mathcal{P}_{p}\left(\mathbb{R}^{d}\right)$ is a Polish space under the $p$-Wasserstein metric.

Regarding the $\cW_p$ distance, we have the following simple observations.
\begin{lmm}\label{lmm:w2aux}
Let $p\ge 1$.
For any probability measure $\mu$,
\[
\mathcal{W}_{p}^{p}(\mu,\delta_{0})
=\int |y|^p \mu(dy).
\]
For two empirical measures
$\mu_1=\frac{1}{N}\sum_{i=1}^N\delta_{Y^{1,i}}$ and $\mu_2=\frac{1}{N}\sum_{i=1}^N \delta_{Y^{2,i}}$. Then,
\[
\mathcal{W}_p^p(\mu_1, \mu_2)
\le \frac{1}{N}\sum_{i=1}^N |Y^{1,i}-Y^{2,i}|^p.
\]
\end{lmm}
The proof for the first claim with $p=2$ can be found, for example, in \cite[Lemma 2.3]{dos2019}, and the proof for general $p\ge 1$ can be similarly proved. The second follows from constructing a simple transport plan $\pi(dx, dy)=\frac{1}{N}\sum_{i=1}^N \delta_{Y^{1i}}(dx)\otimes \delta_{Y^{2i}}(dy)$.

We will make the following technical assumptions on the coefficients, which are essentially the continuity and growth conditions. %\tcr{Borrow from \cite{stockinger2021}.}

\begin{assumption}\label{assp3.3}
There exists a constant $\gamma \ge 0$ such that
\begin{equation}\label{Polynomial growth}
|a(x,\mu)-a(\bar{x},\bar{\mu})|\vee\left\|b(x,\mu)-b(\bar{x},\bar{\mu})\right\|\leq C[(1+|x|^{\gamma}+|\bar{x}|^{\gamma})|x-\bar{x}|+\mathcal{W}_{2}(\mu, \bar{\mu})],
\end{equation}
 for all $x,\bar{x} \in \mathbb{R}^{d}$ and $\mu,\bar{\mu} \in\mathcal{P}_{2}(\mathbb{R}^{d})$.
\end{assumption}

\begin{assumption}\label{assp3.2}
There exists a pair of constants $ L > 0$ and $p_0> \max(2\gamma, 4)$ such that the coefficients $a(\cdot, \cdot): \R^d\times \mathcal{P}_{2}(\R^d)\to \R^d$ and $b(\cdot, \cdot): \R^d\times \mathcal{P}_{2}(\R^d)\to \R^{d\times m'}$  satisfy the initial condition  $|a(0, \mu)|+|b(0, \mu)|\le L\mathcal{W}_2(\mu, \delta_0)$ and
\begin{equation}\label{Khasminskii2}
\begin{aligned}
 (x-\bar{x})\cdot\left(a(x, \mu)-a(\bar{x}, \bar{\mu})\right) +\frac{1}{2}(p_0-1)\left\|b(x, \mu)-b(\bar{x}, \bar{\mu})\right\|^{2}
 \leq L\left\{|x-\bar{x}|^{2}+\mathcal{W}_{2}^{2}(\mu, \bar{\mu})\right\},
\end{aligned}
 \end{equation}
for all x, $\bar{x}\in\mathbb{R}^{d}$ and $\mu$, $  \bar{\mu}\in\mathcal{P}_{2}\left(\mathbb{R}^{d}\right)$.
\end{assumption}

The assumptions here also appeared in \cite{stockinger2021}. We note that the Lipschitz continuity for $\mu$ is uniform. A typical example is $a(x, \mu)=f(x)+\int k(x, y)\mu(dy)$ where $k$ is globally Lipschitz. As another comment,  for condition \eqref{Khasminskii2} to hold, $b$ is not necessarily globally Lipschitz in $x$. In fact, as long as the growth can be controlled by the confining effect in $a$, the condition still holds.
One example is $a=x-x^3$ and $b=x^2$ for $x\in \R$.

The lower bound $2\gamma$ for $p_0$ in Assumption \ref{assp3.2} is essential in our proof and is used in the proof of Theorem \ref{thm3.4} where we need to bound $pq_1\gamma$ moment for some $p\ge 2$ and $q_1>1$. The lower bound $4$ for $p_0$, however, is not essential for the convergence. It only ensures the existence of $p\in [2, p_0/2)$ so the error estimates given below in Proposition \ref{pc} and  Theorem \ref{thm3.4} can be cleaner
(otherwise, there would be a term like $N^{-(p-q)/q}$).

The assumptions made above clearly imply the following growth conditions of the coefficients:
\begin{itemize}
\item
For any $q<p_0$ where $p_0$ is the parameter in Assumption \ref{assp3.2}, there exists a constant $ L > 0$ such that
\begin{equation}\label{Khasminskii1}
x\cdot a(x, \mu)+\frac{1}{2}\left(q-1\right)\left\|b(x, \mu)\right\|^{2} \leq L\left\{(1+|x|)^{2}+\mathcal{W}_{2}^{2}\left(\mu, \delta_{0}\right)\right\} \\
 \end{equation}
for all $x\in\mathbb{R}^{d}$ and $\mu\in\cP_{2}\left(\mathbb{R}^{d}\right)$.

In fact, taking $\bar{x}=0$ and $\bar{\mu}=\delta_0$ in  Assumption \ref{assp3.2}, one has
\[
x\cdot(a(x, \mu)-a(0, \delta_0))+\frac{1}{2}(p_0-1)\left\|b(x, \mu)-b(0, \delta_0)\right\|^2\le \tcb{L}(|x|^2+\mathcal{W}_2^2(\mu, \delta_0)).
\]
Moreover, since $\left\|b(x,\mu)\right\|^2 \le \left\|b(x, \mu)-b(0, \delta_0)\right\|^2+2\left\|b(x, \mu)-b(0, \delta_0)\right\|\left\|b(0, \delta_0)\right\|+\left\|b(0, \delta_0)\right\|^2$
and $2uv\le \delta u^2+\frac{1}{\delta }v^2$ for any $u, v\in \R$ and $\delta>0$, one has
\[
\frac{1}{2}(q-1)\left\|b(x,\mu)\right\|^2
\le \frac{1}{2}(p_0-1)\left\|b(x,\mu)-b(0,\delta_0)\right\|^2+\frac{(q-1)(p_0-1)}{2(p_0-q)}\left\|b(0, \delta_0)\right\|^2.
\]
Combining these two simple estimates, the conclusion then follows.

\item
For the same constant $\gamma$ as in Assumption \ref{assp3.3}, it holds that
\begin{equation}\label{Polynomial growth2}
|a(x, \mu)|\vee \|b(x,\mu)\|\leq C[(1+|x|^{\gamma+1}|)+\mathcal{W}_{2}(\mu,\delta_{0})],
\end{equation}
 for all  $x \in \mathbb{R}^{d}$ and $\mu\in\mathcal{P}_{2}(\mathbb{R}^{d})$.
\end{itemize}

%The proof of this proposition  can be found in \cite{stockinger2021}.

Corresponding to \eqref{IN}, we consider
$\bar{X}_t^i$ which solves \eqref{McKeanLimit} but with $\bar{X}_0^i=X_0^{i}$ and the same Brownian motions. In other words, $\bar{X}_t^i$'s satisfy the following
\begin{equation}\label{I}
\bar{X}_{t}^{i}=X_{0}^{i}+\int_{0}^{t} a\left(\bar{X}_{s}^{i}, \mathcal{L}\left(\bar{X}_{s}^{i}\right)\right) d s+\int_{0}^{t}b\left(\bar{X}_{s}^{i}, \mathcal{L}\left(\bar{X}_{s}^{i}\right)\right) d W_{s}^{i}
 \end{equation}
almost surely for any $t \in[0, T]$ and $i\in{1,\cdots , N}$. Here, $W_t^i$ is the same process as in \eqref{IN} and $X^{i}_{0} \in \R^d$ are the initial values of the SDEs. Clearly, $\{\bar{X}_t^i\}$ are i.i.d. with the law being $\cL(X_t)$ for the McKean-Vlasov SDE \eqref{McKeanLimit}.

For the convenience of the discussion, we denote the moments of $\mu$
for $q\ge 1$ by
\begin{equation*}
M_{q}(\mu):=\int_{\mathbb{R}^{d}}|x|^{q} \mu(d x).
\end{equation*}
Hence, $M_q(\mu_t)=\E[|\bar{X}_t^i|^q]$ for any $i$.

By standard techniques, one may show the following bounds on the moments (see, for example, \cite[Theorem 1]{dereich2013} ). For the convenience of the readers, we provide a short sketch in Appendix \ref{app:bound1}.
\begin{prpstn}\label{pro:momentbound1}
Let \eqref{Khasminskii1} hold for some $q\ge 2$ and $\mu_0\in \cP_q$. Then for any $T > 0$, there exists $C$ depending on $T$ and $q$ but independent of $N$ and $i$ such that
\begin{equation}\label{eq:momentbound1}
\sup _{ 0 \leq t \leq T}\mathbb{E}[|X_{t}^{i}|^{q}+|\bar{X}_t^i|^q]\leq C\left(1+\mathbb{E}[|X_{0}^{i}|^{q}]\right).
\end{equation}
\end{prpstn}

The following result claims that the trajectories of $\bar{X}_t^i$ and $X_t^{i}$ become identical as $N\to\infty$ almost surely. This implies the propagation of chaos for the interacting particle system so that the mean field limit to the McKean SDE holds.

\begin{prpstn}[Propagation of Chaos]\label{pc}
 Let Assumption \ref{assp3.2} be satisfied. If for some $p\in [2, p_0/2)$ the initial law has finite $p$-moment or $\mu_0\in \cP_p \subset\cP_2$, then it holds that
\begin{equation}
\sup _{i \in\{1, \ldots, N\}} \sup _{t \in[0, T]} \E[\left|X_{t}^{i}-\bar{X}_{t}^{i}\right|^{p}] \leq C
\left\{\begin{array}{ll}
N^{-1 / 2}, & \text { if } p>d/2, \\
N^{-1 / 2} \log(1+N), & \text { if } p=d/2,\\
N^{-p / d}, & \text { if }p\in[2,d/2),
\end{array}\right.
\end{equation}
where the constant $C> 0$ depends on $T$, $p$ and $p_0$ but does not depend on N.
\end{prpstn}
\begin{proof}
From equations (\ref{IN}), (\ref{I}), we have
\begin{equation}
\bar{X}_{t}^{i}-X_{t}^{i}=\int_{0}^{t}a\left(\bar{X}_{s}^{i}, \mathcal{L}\left(\bar{X}_{s}^{i}\right)\right)-a\left(X_{s}^{i}, \mu_{s}^{X}\right)ds+\int_{0}^{t}b\left(\bar{X}_{s}^{i}, \mathcal{L}\left(\bar{X}_{s}^{i}\right)\right)-b\left(X_{s}^{i}, \mu_{s}^{X}\right)dW_{s}^{i}.
\end{equation}

Using It\^o’s formula,
\begin{equation*}
\begin{aligned}
\mathbb{E}[|\bar{X}_{t}^{i}-X_{t}^{i}|^{p}]& \leq  ~p \mathbb{E}\left[\int_{0}^{t}|\bar{X}_{s}^{i}-X_{s}^{i}|^{p-2}\left(\bar{X}_{s}^{i}-X_{s}^{i}\right)\cdot \left(a\left(\bar{X}_{s}^{i}, \mathcal{L}\left(\bar{X}_{s}^{i}\right)\right)-a\left(X_{s}^{i}, \mu_{s}^{X}\right)\right) ds\right]\\ &\quad+\frac{p(p-1)}{2} \mathbb{E}\left[\int_{0}^{t}|\bar{X}_{s}^{i}-X_{s}^{i}|^{p-2}\|b\left( \bar{X}_{s}^{i}, \mathcal{L}\left(\bar{X}_{s}^{i}\right)\right)-b\left(X_{s}^{i}, \mu_{s}^{X}\right)\|^2 ds\right].
\end{aligned}
\end{equation*}
Due to Assumption \ref{assp3.2}, one has
\begin{equation*}
\begin{aligned}
\mathbb{E}[|\bar{X}_{t}^{i}-X_{t}^{i}|]^{p}
 &\leq  L \mathbb{E}\left[\int_{0}^{t} |\left(\bar{X}_{s}^{i}-X_{s}^{i}\right)|^{p-2}[\mathcal{W}_{2}^{2}\left(\mathcal{L}\left(\bar{X}_{s}^{i}\right), \mu_{s}^{X}\right) +|\bar{X}_{s}^{i}-X_{s}^{i}|^{2}]ds\right] \\
 & \leq  L \mathbb{E} \left[\int_{0}^{t}|\left(\bar{X}_{s}^{i}-X_{s}^{i}\right)|^{p}ds\right]+L\mathbb{E} \left[\int_{0}^{t}|\left(\bar{X}_{s}^{i}-X_{s}^{i}\right)|^{p-2}\mathcal{W}_{2}^{2}\left(\mathcal{L}\left(\bar{X}_{s}^{i}\right), \mu_{s}^{X}\right)ds\right]\\
& \leq C\mathbb{E} \left[\int_{0}^{t}\left|\left(\bar{X}_{s}^{i}-X_{s}^{i}\right)\right|^{p}ds\right]+C \mathbb{E}\left[ \int_{0}^{t}
\mathcal{W}_{2}^{p}\left(\mathcal{L}\left(\bar{X}_{s}^{i}\right), \mu_{s}^{X}\right)ds\right],\\
\end{aligned}
\end{equation*}
for any $t \in [0, T]$.  The last inequality is due to the Young's inequality.

Denote $\mu_s^{\bar{X}}:=  \frac{1}{N} \sum_{j=1}^{N} \delta_{\bar{X}_{t}^{j}}$.  By Lemma \ref{lmm:w2aux}, one has
\begin{equation*}
\begin{aligned}
\mathcal{W}_{2}^{p}\left(\mathcal{L}\left(\bar{X}_{s}^{i}\right), \mu_{s}^{X}\right)&=\left[\mathcal{W}_{2}^{2}\left(\mathcal{L}\left(\bar{X}_{s}^{i}\right), \mu_{s}^{X}\right)\right]^{p/2} \\
&\le  \left[2\mathcal{W}_2(\mu_s^{\bar{X}}, \mu_s^X)^2+2\mathcal{W}_2^2(\mu_s^{\bar{X}}, \cL(\bar{X}_s^1)) \right]^{p/2} \\
&\leq C \left[\frac{1}{N} \sum_{j=1}^{N}\left|\bar{X}_{s}^{j}-X_{s}^{j}\right|^{2}\right]^{p/2}+C \mathcal{W}_{2}^{p}\left(\frac{1}{N} \sum_{j=1}^{N} \delta_{\bar{X}_{t}^{j}}, \mathcal{L}\left(\bar{X}_{s}^{1}\right)\right).
\end{aligned}
\end{equation*}

Consequently, one obtains
 \begin{equation*}
\begin{aligned}
\mathbb{E}[\left|\bar{X}_{t}^{i}-X_{t}^{i}\right|^{p}]&\leq C \int_{0}^{t}\mathbb{E}[\left|\left(\bar{X}_{s}^{i}-X_{s}^{i}\right)\right|^{p}]ds+C \int_{0}^{t}\mathbb{E}\left[\frac{1}{N} \sum_{j=1}^{N}\left|\bar{X}_{s}^{j}-X_{s}^{j}\right|^{2}\right]^{p/2}ds
\\&\quad+C\int_{0}^{t}\mathbb{E}\left[
\mathcal{W}_{2}^{p}\left(\frac{1}{N} \sum_{j=1}^{N} \delta_{\bar{X}_{t}^{j}}, \mathcal{L}\left(\bar{X}_{s}^{1}\right)\right)\right]ds\\
&\leq C \int_{0}^{t}\mathbb{E}[\left|\bar{X}_{s}^{i}-X_{s}^{i}\right|^{p}]ds+C\int_{0}^{t}\mathbb{E}\left[
\mathcal{W}_{p}^{p}\left(\frac{1}{N} \sum_{j=1}^{N} \delta_{\bar{X}_{t}^{j}}, \mathcal{L}\left(\bar{X}_{s}^{1}\right)\right)\right]ds.
\end{aligned}
\end{equation*}
In the last inequality, we have applied  $\mathcal{W}_2(\mu, \nu)
\le \mathcal{W}_p(\mu, \nu)$. Then, by  the Minkowski inequality,
\[
\mathbb{E}\left[\frac{1}{N} \sum_{j=1}^{N}\left|\bar{X}_{s}^{j}-X_{s}^{j}\right|^{2}\right]^{p/2}
\le \left[\frac{1}{N}\sum_{j=1}^N \|\left|\bar{X}_{s}^{j}-X_{s}^{j}\right|^{2}\|_{L^{p/2}} \right]^{p/2}
=\left[\frac{1}{N}\sum_{j=1}^N (\E[\left|\bar{X}_{s}^{j}-X_{s}^{j}\right|^{p}])^{2/p} \right]^{p/2}
=\E[\left|\bar{X}_{s}^{i}-X_{s}^{i}\right|^{p}].
\]
The last equality here holds due to the symmetry.

The term $\mathbb{E}[\mathcal{W}_{p}^{p}(\frac{1}{N} \sum_{j=1}^{N} \delta_{\bar{X}_{t}^{j}}, \mathcal{L}\left(\bar{X}_{s}^{1}\right))]$ is controlled by the Wasserstein distance estimate for the empirical measures with i.i.d. samples  in  \cite[Theorem 1]{fournier2015}.
Here, since $p_0>2p$, one can choose the moment $q$ such that $q>2p$ and the terms $N^{-(p-q)/q}$ in  \cite[Theorem 1]{fournier2015} can be removed for simplicity. Then, applying Gr\"onwall’s inequality completes the proof.
\end{proof}

We remark that the requirement $p<p_0/2$ is only used such that the term $N^{-(p-q)/q}$ is absent for simplicity. Clearly, if $p\in [p_0/2, p_0)$, there is still convergence.
Another remark is that the rate of the propagation of chaos here is like $N^{-1/(2p)}$ when $p$ is large. This low rate is due to the usage of Wasserstein distances to gauge the continuity on the distribution dependence in the coefficients.
If the dependence on the distribution is through some statistics, the rate may be improved to $N^{-1/2}$ due to the central limit theorem. Note that though the particles are dependent for finite $N$, the rate of central limit theorem is expected to hold in the large $N$ regime due to propagation of chaos, though rigorous proof needs careful estimates (see, for example, \cite{wang2021gaussian} regarding the rate of convergence for statistics).

\section{A truncated Euler scheme}\label{part3}

The drift and diffusion coefficients can be unbounded. To overcome this difficulty, we make use of the truncation techniques \cite{guo2018,mao2015,mao2016convergence}.
For this purpose, we choose a strictly increasing continuous function $\varphi$ such
that
\[
\sup_{|x| \leq u}\left(|a(x)| \vee\left\|b(x)\right\| \right) \leq \varphi(u).
\]
The inverse function $\varphi^{-1}$ is strictly increasing continuous function from $[\varphi(0),+\infty)$
to $\mathbb{R}_{+}$. Firstly, define a discretization grid $T_ {N_{T}}: 0 = t_ {0} < t_ {1} < \cdots < t_ {N_{T}} =T $, where $t_ {m} =m \Delta $ with the step size $ \Delta=T/N_{T} $ and $N_{T} \in \mathbb {N}$ is the number of subintervals. Take $\Delta^{*}\in (0,1]$, which is meant to be the largest step size, and a strictly decreasing function $h : (0,\Delta^{*}]\longrightarrow(0,+\infty)$ such that
 \begin{equation*}
h\left(\Delta^{*}\right) \geq \varphi(1), \quad \lim _{\Delta \rightarrow 0} h(\Delta)=\infty \quad\text { and }\quad  \sup_{ \Delta \in\left(0, \Delta^{*}\right]}\Delta^{1/4} h(\Delta) <\infty.
 \end{equation*}
One possible example is $h(\Delta)=\Delta^{-\frac{\varepsilon}{2}}$ and $\epsilon$ is any value in $(0, 1/4]$. Another possible example is $h(\Delta)=|\log \Delta|$.
 For a given step size $\Delta\in (0,1]$, define the truncated functions by
 \begin{equation*}
 a^{\Delta}(x,\mu)=a\left(\left(|x| \wedge \varphi^{-1}(h(\Delta))\right) \frac{x}{|x|}, \quad \mu\right),\quad
 b^{\Delta}(x,\mu)=b\left(\left(|x| \wedge \varphi^{-1}(h(\Delta))\right) \frac{x}{|x|}, \quad \mu \right).
\end{equation*}

Similar to that in \cite{mao2015}, the truncated coefficients satisfy the following properties.
\begin{lmm}
It holds that
\begin{equation}\label{Hcondition}
|a^{\Delta}(x,\mu)|\vee\|b^{\Delta}(x,\mu)\|\leq\varphi(\varphi^{-1}(h(\Delta)))=h(\Delta),
\end{equation}
and for any $q\in [2, p_0)$ where $p_0$ is  parameter in Assumption \ref{assp3.2}, there exists a constant $ L > 0$ such that
\begin{equation}\label{DeltaKha}
x\cdot a^{\Delta}(x, \mu)+\frac{1}{2}\left(q-1\right)\left\|b^{\Delta}(x, \mu)\right\|^{2} \leq L\left\{(1+|x|)^{2}+\mathcal{W}_{2}^{2}\left(\mu, \delta_{0}\right)\right\} ,
\quad \forall x\in\mathbb{R}^d, \mu\in\mathcal{P}_{2}\left(\mathbb{R}^d\right).
\end{equation}
\end{lmm}

\begin{proof}
Here, we only need to verify the Khasminskii-type condition \eqref{DeltaKha}.

 When $|x|\le \varphi^{-1}(h(\Delta))$, it is simply \eqref{Khasminskii1}.
For $|x|> \varphi^{-1}(h(\Delta))$,  due to the dependence on $\mu$, unlike the argument in \cite{mao2015}, we start by \eqref{Khasminskii2}  with
\[
x\cdot a(x, \mu)+\frac{1}{2}(p_0-1) \left\|b(x, \mu)-b(0, \mu)\right\|^{2}  \le x\cdot a(0, \mu)+L|x|^2.
\]

Denote
\[
\bar{x}=\varphi^{-1}(h(\Delta))\frac{x}{|x|}.
\]
Then, it follows that
\begin{equation*}
\begin{split}
x\cdot a^{\Delta}(x, \mu)+\frac{1}{2}(p_0-1) \left\|b^{\Delta}(x, \mu)-b^{\Delta}(0, \mu)\right\|^{2}
&=x\cdot a(\bar{x}, \mu)+\frac{1}{2}(p_0-1) \left\|b(\bar{x}, \mu)-b(0, \mu)\right\|^{2}\\
&\le \frac{|x|}{\varphi^{-1}(h(\Delta))}
\left[\bar{x}\cdot a(\bar{x}, \mu)+\frac{1}{2}(p_0-1)\left\|b(\bar{x}, \mu)-b(0, \mu)\right\|^{2} \right] \\
& \le \frac{|x|}{\varphi^{-1}(h(\Delta))}
\left( \bar{x}\cdot a(0, \mu)+L(\varphi^{-1}(h(\Delta)))^2\right) \\
& =x\cdot a(0, \mu)+L|x|\varphi^{-1}(h(\Delta))\\
& \le x\cdot a(0, \mu)+L|x|^2.
\end{split}
\end{equation*}

%Further, \tcr{given the initial condition  $|a(0, \mu)|+|b(0, \mu)|\le C\mathcal{W}_2(\mu, \delta_0)$}.
Then a similar argument as in deriving \eqref{Khasminskii1} gives the desired result.
\end{proof}

\subsection{The truncated scheme}

The numerical approximation of $X_{t}^{i}$ at $t_m$ is denoted by $X_{t_{m}}^{i,\Delta}$.
With the truncated coefficients, the numerical solutions $X_{t_{m}}^{i,\Delta}$ are then generated by the basic Euler-Maruyama scheme
\begin{equation}\label{eq:basicscheme}
 X_{t_{m+1}}^{i,\Delta}=X_{t_{m}}^{i,\Delta}+ a^{\Delta}\left(X_{t_{m}}^{i,\Delta}, \mu_{t_{m}}^{X,\Delta}\right)\Delta + b^{\Delta}\left(X_{t_{m}}^{i,\Delta}, \mu_{t_{m}}^{X,\Delta}\right) \Delta W_{t_{m}}^{i},
\end{equation}
where $\Delta W_{t_{m}}^{i}=W_{t_{m+1}}^{i}-W_{t_{m}}^{i}$, the initial value $X_{0}^{i,\Delta}=X^{i}_{0}$ and
\[
\mu_{t_{m}}^{X, \Delta}(\cdot):=\frac{1}{N}\sum_i \delta_{X_{t_{m}}^{i,\Delta}}(\cdot).
\]
The scheme \eqref{eq:basicscheme} is the basic truncated scheme studied in this work.

We introduce two versions of extension of the numerical solution at the discrete time points to $t\in [0,\infty)$. The first is the piecewise constant extension given by
\begin{equation*}
\hat{X_{t}}^{i,\Delta}= X_{t_{m}}^{i,\Delta}, \quad t_{m}\leq t < t_{m+1},
\end{equation*}
and the second is the continuous extension of the truncated Euler-Maruyama method defined by
\begin{equation}\label{INn}
 X_{t}^{i, \Delta}=X_{t_{m}}^{i,\Delta}+\int_{t_{m}}^{t} a^{\Delta}\left(\hat{X}_{s}^{i,\Delta}, \hat{\mu}_{s}^{X,\Delta}\right) d s+\int_{t_{m}}^{t} b^{\Delta}\left(\hat{X}_{s}^{i,\Delta}, \hat{\mu}_{s}^{X,\Delta}\right) d W_{s}^{i}.
\end{equation}
Here, the initial value $\hat{X}_{0}^{i,\Delta}=X^{i}_{0}$ and
\[
\hat{\mu}_t^{X,\Delta}(\cdot):=\frac{1}{N}\sum_{i=1}^N \delta_{\hat{X}_t^{i,\Delta}}(\cdot).
\]

\subsection{Analysis of the truncated scheme}

We present some preliminary properties of the numerical solutions. We first give some moments estimates of the numerical solution to ensure certain stability.  The first lemma below is to control the difference between the two versions of extension.
\begin{lmm}\label{lmm:auxiliary1}
For any $\Delta\in (0,\Delta^{*}]$, any $t \geq 0$ and any $p \geq 2$,
\begin{equation}
\mathbb{E}[|X_{t}^{i,\Delta}-\hat{X}_{t}^{i,\Delta}|^{p}] \leq C \Delta^{p/2}h^{p}(\Delta),
\end{equation}
where C is a positive constant independent of $\Delta$.
\end{lmm}
\begin{proof}
By the definitions \eqref{INn} , one has
\[
X_{t}^{i,\Delta}-\hat{X}_{t}^{i,\Delta}
=\int_{t_{m}}^{t} a^{\Delta}\left(\hat{X}_{s}^{i,\Delta}, \hat{\mu}_{s}^{X,\Delta}\right) d s+\int_{t_{m}}^{t} b^{\Delta}\left(\hat{X}_{s}^{i,\Delta}, \hat{\mu}_{s}^{X,\Delta}\right) d W_{s}^{i}.
\]

Using H\"older’s inequality, Burkholder-Davis-Gundy inequality and (\ref{Hcondition}), we have
\begin{equation*}
\begin{aligned}
 \mathbb{E}[|X_{t}^{i,\Delta}-\hat{X}_{t}^{i,\Delta}|^{p}] &\leq  C \mathbb{E}\left[\left|\int_{t_{m}}^{t} a^{\Delta}\left(\hat{X}_{s}^{i,\Delta}, \hat{\mu}_{s}^{X,\Delta}\right) d s\right|^{p}+\left|\int_{t_{m}}^{t} b^{\Delta}\left(\hat{X}_{s}^{i,\Delta}, \hat{\mu}_{s}^{X,\Delta}\right) d W_{s}^{i}\right|^{p}\right]\\
 &\leq  C\Delta^{p-1} \mathbb{E}\left[\int_{t_{m}}^{t} \left|a^{\Delta}\left(\hat{X}_{s}^{i,\Delta},\hat{\mu}_{s}^{X,\Delta}\right)\right|^{p} d s\right]+C\Delta^{(p-2)/2}\mathbb{E}\left[\int_{t_{m}}^{t} \left\|b^{\Delta}\left(\hat{X}_{s}^{i,\Delta}, \hat{\mu}_{s}^{X,\Delta}\right)\right\|^{p} d s\right]\\
 &\leq  C\Delta^{p-1} \mathbb{E}\left[\int_{t_{m}}^{t} h^{p}(\Delta) d s\right]+C\Delta^{(p-2)/2}\mathbb{E}\left[\int_{t_{m}}^{t} h^{p}(\Delta) d s\right]\\
 & \leq  C\Delta^{p/2} h^{p}(\Delta).
 \end{aligned}
 \end{equation*}
\end{proof}

\begin{prpstn}\label{prpstn:boundedness1}
Suppose \eqref{DeltaKha} holds and $\mu_0\in \cP_q \subset \cP_2$. Then for any $\Delta\in (0,\Delta^{*}]$, any $T > 0$, it holds that
 \begin{equation}
%\tcr{\sup _{0<\Delta \leq \Delta^{*}}}
\sup _{ 0 \leq t \leq T}\mathbb{E}[|X_{t}^{i,\Delta}|^{q}] \leq C\left(1+\mathbb{E}[|X_{0}^{i}|^{q}]\right),
\end{equation}
where C is a positive constant dependent on T and $q$ but independent of $\Delta$.
\end{prpstn}

\begin{proof}
Recall \eqref{INn}:
 \begin{equation*}
d X_{t}^{i, \Delta}= a^{\Delta}\left(\hat{X}_{t}^{i,\Delta}, \hat{\mu}_{t}^{X,\Delta}\right) dt+ b^{\Delta}\left(\hat{X}_{t}^{i,\Delta}, \hat{\mu}_{t}^{X,\Delta}\right) d W_{t}^{i}.
\end{equation*}
By the It\^o's formula, one has
\begin{equation*}
\frac{d}{dt}\mathbb{E}[|X_{t}^{i,\Delta}|^{q}]
 \leq   q \mathbb{E}\left[ |X_{t}^{i,\Delta}|^{q-2}\left( X_{t}^{i,\Delta}\cdot  a^{\Delta}\left(\hat{X}_{t}^{i,\Delta}, \hat{\mu}_{t}^{X,\Delta}\right)\right)\right]
+q\mathbb{E}\left[\frac{q-1}{2}|X_{t}^{i,\Delta}|^{q-2}\left\| b^{\Delta}\left(\hat{X}_{t}^{i,\Delta}, \hat{\mu}_{t}^{X,\Delta}\right)\right\|^{2}\right].
\end{equation*}
Writing $X_{t}^{i,\Delta}=(X_{t}^{i,\Delta}-\hat{X}_t^{i,\Delta})+\hat{X}_t^{i,\Delta}$, one has according to \eqref{DeltaKha} that
\begin{multline*}
\left( X_{t}^{i,\Delta}\cdot  a^{\Delta}\left(\hat{X}_{t}^{i,\Delta}, \hat{\mu}_{t}^{X,\Delta}\right)\right)
+\frac{q-1}{2}\left\| b^{\Delta}\left(\hat{X}_{t}^{i,\Delta}, \hat{\mu}_{t}^{X,\Delta}\right)\right\|^{2}\\
  \le \left( X_{t}^{i,\Delta}-\hat{X}_t^{i,\Delta}\right)\cdot  a^{\Delta}\left(\hat{X}_{t}^{i,\Delta}, \hat{\mu}_{t}^{X,\Delta}\right)
+L \left[(1+|\hat{X}_{t}^{i,\Delta}|)^{2}+\mathcal{W}_{2}^{2}(\hat{\mu}_{t}^{X,\Delta},\delta_{0}) \right].
\end{multline*}
Hence,
\[
\begin{aligned}
\frac{d}{dt}\mathbb{E}[|X_{t}^{i,\Delta}|^{q}] &\le C\mathbb{E}\left[|X_{t}^{i,\Delta}|^{q-2}((1+|\hat{X}_{t}^{i,\Delta}|)^{2}+\mathcal{W}_{2}^{2}(\hat{\mu}_{t}^{X,\Delta},\delta_{0}))\right]
+\E\left[|X_{t}^{i,\Delta}-\hat{X}_{t}^{i,\Delta}|^{q/2}\left| a^{\Delta}\left(\hat{X}_{t}^{i,\Delta}, \hat{\mu}_{t}^{X,\Delta}\right)\right|^{q/2}\right]\\
&\le C(1+\mathbb{E}[|X_{t}^{i,\Delta}|^{q}])
+C\E\left[|X_{t}^{i,\Delta}-\hat{X}_{t}^{i,\Delta}|^{q/2}| a^{\Delta}\left(\hat{X}_{t}^{i,\Delta},\hat{\mu}_{t}^{X,\Delta}\right)|^{q/2}\right]
+\mathbb{E}\left[\left(\frac{1}{N} \sum_{j=1}^{N}|\hat{X}_{t}^{j,\Delta}|^{2}\right)^{q/2}\right],
\end{aligned}
\]
where we applied Lemma \ref{lmm:w2aux}.

\noindent By  Lemma \ref{lmm:auxiliary1}, it is clear that
\[
\E\left[|X_{t}^{i,\Delta}-\hat{X}_{t}^{i,\Delta}|^{q/2}| a^{\Delta}\left(\hat{X}_{t}^{i,\Delta},\hat{\mu}_{t}^{X,\Delta}\right)|^{q/2}\right]
\lesssim h^{q/2}\E\left[|X_{t}^{i,\Delta}-\hat{X}_{t}^{i,\Delta}|^{q/2}\right]
\lesssim \Delta^{q/4}h^{q}(\Delta).
\]

\noindent By the Minkowski inequality and the symmetry
\[
\mathbb{E}\left(\frac{1}{N} \sum_{j=1}^{N}|\hat{X}_{t}^{i,\Delta}|^{2}\right)^{q/2}
\le \left(\frac{1}{N} \sum_{j=1}^{N}(\mathbb{E}[|\hat{X}_{t}^{j,\Delta}|^{q}])^{2/q}\right)^{q/2}=\mathbb{E}[|\hat{X}_{t}^{i,\Delta}|^{q}],\quad \forall i.
\]

\noindent It follows from the fact $|\hat{X}_t| \le \sup_{s\le t} |X_{s}^{i,\Delta}|$ by definition that
\[
\mathbb{E}[|X_{t}^{i,\Delta}|^{q}]\leq (\mathbb{E}[|X_{0}^{i}|^{q}]+C_{T})+C\int_{0}^{t}\sup _{0\leq u\leq s}\mathbb{E}[|X_{u}^{i,\Delta}|^{q}]ds.
\]

\noindent
By the monotonticy of the bound on the right hand side, one has
\begin{equation*}
\sup _{0\leq s\leq t}\mathbb{E}[|X_{s}^{i,\Delta}|^{q}]\leq (\mathbb{E}[|X_{0}^{i}|^{q}]+C_{T})+C\int_{0}^{t}\sup _{0\leq u\leq s}\mathbb{E}[|X_{u}^{i,\Delta}|^{q}]ds.
\end{equation*}
Application of Gr\"onwall’s inequality then completes the proof.
\end{proof}

Below, we prove the convergence of the truncated scheme. Define two stopping times for particle $i$
\begin{equation*}
\tau_{i}:=\inf \{t \geq 0,|\bar{X}_{t}^{i}| \geq \varphi^{-1}(h(\Delta))\}
\end{equation*} and
\begin{equation*}
\rho_{i}:=\inf \{t \geq 0,|X_{t}^{i,\Delta}| \geq \varphi^{-1}(h(\Delta))\}.
\end{equation*}
 Let $\theta_{i}=\tau_{i}\wedge\rho_{i}$.  The estimates in Propostion \ref{prpstn:boundedness1} ensures that $\theta_{i}\to +\infty$ a.s. as $\Delta\to 0$. Clearly, before $\theta_{i}$, the truncated method reduces to the usual Euler–Maruyama scheme.

 Below, we make use of these stopping times to obtain the convergence of the truncated scheme, and the following is our first main result.

\begin{thrm}[Rate of Convergence]\label{thm3.4}
Suppose Assumption \ref{assp3.3} and Assumption \ref{assp3.2} hold and $\mu_0\in \cP_{p_{0}}$. Let $\{X_t^i\}$
be the solution to the $N$-particle system \eqref{IN}  and $X_t^{i, \Delta}$ be the solution to the truncated scheme \eqref{INn}. Denote $\eta(\Delta):=1/(\varphi^{-1}\circ h(\Delta))$.  Then, for $p\in [2, p_0/2)$ with $p\gamma <p_0$, it holds for every $q<p_0$ that
\begin{equation}\label{eq:convest1}
 \sup_{0\leq t\leq T}\mathbb{E}[| X_{t}^{i}-X_{t}^{i,\Delta}|^{p}]\leq C_{T}\Delta^{p/2}h^p(\Delta)+C_{T}\eta(\Delta)^{q-p},
\end{equation}
where $C_{T}$ is a positive constant dependent on $p, q, T$ but independent of $N$ and $\Delta$. Consequently,
one has for $q<p_0$,
\begin{equation}\label{eq:convest2}
 \sup_{t\le T}\mathbb{E}[|\bar{X}_{t}^{i}-X_{t}^{i,\Delta}|^{p}] \leq
 C\left\{\begin{array}{ll}N^{-1 / 2}+\Delta^{p/2}h(\Delta)^{p}+
 \eta(\Delta)^{q-p}, & \text { if } p>d/2, \\
  N^{-1 / 2} \log(1+N)+\Delta^{p/2}h(\Delta)^{p}+\eta(\Delta)^{q-p}, & \text { if } p=d/2, \\
 N^{-p/ d}+\Delta^{p/2} h(\Delta)^{p}+\eta(\Delta)^{q-p}, & \text { if }p\in[2,d/2),
 \end{array}\right.
\end{equation}
where $C$ depends on $p,q,T$ but is independent of $N$ and $\Delta$. Hence, it holds that
\[
 \lim_{N\to\infty, \Delta\to 0}\sup_{t\le T}\left(\mathbb{E}[|\bar{X}_{t}^{i}-X_{t}^{i,\Delta}|^{p}]\right)^{1/p}=0.
\]
\end{thrm}
\begin{proof}
For $0\leq s\leq t\wedge \theta_{i}$, the coefficients $a$ and $b$ that particle $i$ sees  agree with the truncated ones, $a^{\Delta}$ and $b^{\Delta}$. Hence,  by the It\^o's formula,
 \begin{equation*}
 \begin{aligned}
\mathbb{E}[|X_{t\wedge\theta_{i}}^{i}-X_{t\wedge\theta_{i}}^{i,\Delta}|^{p}]
&\le  p\mathbb{E}\left[\int_{0}^{t\wedge\theta_{i}}|X_{s}^{i}-X_{s}^{i,\Delta}|^{p-2}\left( X_{s}^{i}
-X_{s}^{i,\Delta}\right)\cdot\left( a(X_{s}^{i},\mu_{s}^{X})-a(\hat{X}_{s}^{i,\Delta},\hat{\mu}_{s}^{X,\Delta})\right) ds\right]
\\&\quad+\frac{p(p-1)}{2}\mathbb{E}\left[\int_{0}^{t\wedge\theta_{i}}|X_{s}^{i}-X_{s}^{i,\Delta}|^{p-2}\|b( X_{s}^{i},\mu_{s}^{X})-b(\hat{X}_{s}^{i,\Delta},\hat{\mu}_{s}^{X,\Delta})\|^{2}ds\right],
\end{aligned}
\end{equation*}

By simple splitting, $$\left( X_{s}^{i}
-X_{s}^{i,\Delta}\right)\cdot\left(a(X_{s}^{i},\mu_{s}^{X})-a(\hat{X}_{s}^{i,\Delta},\hat{\mu}_{s}^{X,\Delta})\right)+\frac{p-1}{2}\|b(X_{s}^{i},\mu_{s}^{X})-b(\hat{X}_{s}^{i,\Delta}, \hat{\mu}_{s}^{X,\Delta})\|^{2}\le A_1+A_2,$$
with
\[
A_1 =\left( X_{s}^{i}
-X_{s}^{i,\Delta}\right)\cdot\left(a(X_{s}^{i},\mu_{s}^{X})-a(X_{s}^{i,\Delta},\mu_{s}^{X,\Delta})\right)+ \frac{p'-1}{2}\|b(X_{s}^{i},\mu_{s}^{X})-b(X_{s}^{i,\Delta},\mu_{s}^{X,\Delta})\|^{2}
\]
and
\[
A_2=\left( X_{s}^{i}
-X_{s}^{i,\Delta}\right)\cdot\left(a(X_{s}^{i,\Delta},\mu_{s}^{X,\Delta})-a(\hat{X}_{s}^{i,\Delta}, \hat{\mu}_{s}^{X,\Delta})\right)
+C(p',p)\|b(X_{s}^{i,\Delta},\mu_{s}^{X,\Delta})-b(\hat{X}_{s}^{i,\Delta},\hat{\mu}_{s}^{X,\Delta})\|^{2}.
\]
Here, we can take any $p'\in (p, p_0)$.

According to Assumption \ref{assp3.2}.
\[
A_1 \le C\left(| X_{s}^{i}-X_{s}^{i,\Delta}|^{2}
+\mathcal{W}_{2}^{2}(\mu_{s}^{X},\mu_{s}^{X,\Delta})\right).
\]
The terms in $A_2$ can be controlled using Assumption \ref{assp3.3}:
\[
A_2\lesssim
|X_{s}^{i}-X_{s}^{i,\Delta}|^{2}
+(1+|X_{s}^{i,\Delta}|^{\gamma}+|\hat{X}_{s}^{i,\Delta}|^{\gamma})^2|X_{s}^{i,\Delta}-\hat{X}_{s}^{i,\Delta}|^{2}
+\mathcal{W}_{2}^{2}(\mu_{s}^{X,\Delta}, \hat{\mu}_{s}^{X,\Delta})
\]

Applying Young's inequality, one has
\begin{equation*}
\begin{aligned}
\mathbb{E}[|X_{t\wedge\theta_{i}}^{i}-X_{t\wedge\theta_{i}}^{i,\Delta}|^{p}]
&\lesssim  \left\{\mathbb{E}\left[\int_{0}^{t\wedge\theta_{i}}|X_{s}^{i}-X_{s}^{i,\Delta}|^{p}ds\right]
+\mathbb{E}\left[\int_{0}^{t\wedge\theta_{i}}\mathcal{W}_{2}^{p}(\mu_{s}^{X},\mu_{s}^{X,\Delta})ds\right]\right\}
\\
&\quad+
\mathbb{E}\left[\int_{0}^{t\wedge\theta_{i}}(1+|X_{s}^{i,\Delta}|^{\gamma}+|\hat{X}_{s}^{i,\Delta}|^{\gamma})^p|X_{s}^{i,\Delta}-\hat{X}_{s}^{i,\Delta}|^{p}ds\right]
\\
&\quad+\mathbb{E}\left[\int_{0}^{t\wedge\theta_{i}}\mathcal{W}_{2}^{p}(\mu_{s}^{X,\Delta}, \hat{\mu}_{s}^{X,\Delta})ds\right]=: J_1+J_2+J_3.
\end{aligned}
\end{equation*}

By Lemma \ref{lmm:w2aux}, the Minkowski inequality, and Lemma \ref{lmm:auxiliary1}
\[
J_3
\le \E\left[\int_0^{t}
\left(\frac{1}{N}\sum_{j=1}^{N} \left|X_{s}^{j, \Delta}-\hat{X}_{s}^{j, \Delta}\right|^{2}\right)^{p/2}ds\right]
\le \int_0^t \left(\frac{1}{N}\sum_{j=1}^N (\E[|X_{s}^{j, \Delta}-\hat{X}_{s}^{j, \Delta}|^p])^{2/p}\right)^{p/2}ds
\lesssim \Delta^{p/2}h^p(\Delta).
\]

For $J_{2}$, by the H\"older inequality and Lemma \ref{lmm:auxiliary1} and Proposition \ref{prpstn:boundedness1}, one has
\begin{equation}\label{J3}
\begin{aligned}
J_{2}&\leq    \mathbb{E}\left[\int_{0}^{t}(1+|\hat{X}_{s}^{i,\Delta}|^{\gamma}+|X_{s}^{i,\Delta}|^{\gamma})^p|X_{s}^{i,\Delta}-\hat{X}_{s}^{i,\Delta}|^{p}ds\right]\\
& \lesssim
\int_0^t\left(\E \left[(1+|\hat{X}_{s}^{i,\Delta}|^{\gamma}+|X_{s}^{i,\Delta}|^{\gamma})^{pq_1}\right]\right)^{1/q_1}
\left(\E [|X_{s}^{i,\Delta}-\hat{X}_{s}^{i,\Delta}|^{p_1p}]\right)^{1/p_1}\\
& \le C_{T}\Delta^{p/2}h^p(\Delta).
\end{aligned}
\end{equation}
Here, $1/p_1+1/q_1=1$. Note that Lemma \ref{lmm:auxiliary1} holds for any $p_1 p\ge 2$ so we can choose $p_1$ large enough such that
$pq_1\gamma<p_0$. Such $q_1$ exists because $p\gamma<p_0$.

For $J_1$, we first note that
\[
J_1\le
\mathbb{E}\left[\int_{0}^{t}|X_{s\wedge\theta_{i}}^{i}-X_{s\wedge\theta_{i}}^{i,\Delta}|^{p}ds\right]
+\mathbb{E}\left[\int_{0}^{t\wedge\theta_{i}}\mathcal{W}_{2}^{p}(\mu_{s}^{X},\mu_{s}^{X,\Delta})ds\right],
\]
because for $s>\theta_{i}$, the integrand is nonnegative.
Then again, by Lemma \ref{lmm:w2aux}, the Minkowski inequality and the symmetry of the particles,
\[
\begin{split}
\mathbb{E}\left[\int_{0}^{t\wedge\theta_{i}}\mathcal{W}_{2}^{p}(\mu_{s}^{X},\mu_{s}^{X,\Delta})ds\right]
&\leq\mathbb{E}\left[\int_{0}^{t}\mathcal{W}_{2}^{p}(\mu_{s}^{X},\mu_{s}^{X,\Delta})ds\right]\\
&\le \int_0^t \left(\frac{1}{N}\sum_{j=1}^N (\E[|X_{s}^{j}-X_{s}^{j,\Delta}|^p])^{2/p}\right)^{p/2}ds\\
&= \int_{0}^{t}\mathbb{E}[ |X_{s}^{i}-X_{s}^{i,\Delta}|^{p}]ds.
\end{split}
\]
Here, $\theta_i$ has been thrown away because it is not suitable for stopping other particles $X_s^j$.

Therefore,
\begin{align}
\mathbb{E}[|X_{t\wedge\theta_{i}}^{i}-X_{t\wedge\theta_{i}}^{i,\Delta}|^{p}]
&\lesssim \int_{0}^{t}\mathbb{E}[|X_{s\wedge\theta_{i}}^{i}-X_{s\wedge\theta_{i}}^{i,\Delta}|^{p}]ds+\int_{0}^{t}\mathbb{E}[ |X_{s}^{i}-X_{s}^{i,\Delta}|^{p}]ds+\Delta^{p/2}h^p(\Delta).\notag
\end{align}
Note that for all $s\in [0, t]$,
\begin{equation}\label{eq:partition}
\begin{split}
\mathbb{E}[| X_{s}^{i}-X_{s}^{i,\Delta}|^{p}]
&= \mathbb{E}[| X_{s\wedge\theta_{i}}^{i}-X_{s\wedge\theta_{i}}^{i,\Delta}|^{p}\mathcal{I}_{\{\theta_{i}>T\}}]
+\mathbb{E}[| X_{s}^{i}-X_{s}^{i,\Delta}|^{p}\mathcal{I}_{\{\tau_{i}\leq T \text{ or }\rho_{i}\leq T\}}].\\
& \le \mathbb{E}[| X_{s\wedge\theta_{i}}^{i}-X_{s\wedge\theta_{i}}^{i,\Delta}|^{p}]
+\mathbb{E}[| X_{s}^{i}-X_{s}^{i,\Delta}|^{p}\mathcal{I}_{\{\tau_{i}\leq T \text{ or }\rho_{i}\leq T\}}].
\end{split}
\end{equation}
Here, recall that $\mathcal{I}_{G}$ is the indicator function of the set $G$.
By Young's inequality:
\[
\mathbb{E}[| X_{s}^{i}-X_{s}^{i,\Delta}|^{p}\mathcal{I}_{\{\tau_{i}\leq T \text{ or }\rho_{i}\leq T\}}] \le \frac{\delta}{p_2}\E[| X_{s}^{i}-X_{s}^{i,\Delta}|^{pp_2}]+\frac{1}{q_2\delta^{1/(p_2-1)}}\P(\tau_{i}\leq T \text{ or }\rho_{i}\leq T).
\]
Where $p_2>1$ is some positive constant such that $pp_2<p_0$. The expectation in the first term is then bounded.
By the moment estimation of $\bar{X}_{t}^{i}$ and $X_{t}^{i,\Delta}$ and the Markov inequality
\[
\mathbb{P}(\tau_{i}\leq T \text{ or }\rho_{i}\leq T)\leq\frac{\mathbb{E}[|\bar{X}_{t}^{i}|^{q}]+\mathbb{E}[|X_{t}^{i,\Delta}|^{q}]}{(\varphi^{-1}\circ h(\Delta))^{q}},
\]
for $q<p_0$. Denote $\eta(\Delta)=(\varphi^{-1}\circ h(\Delta))^{-1}$. Balancing the terms, we may choose $pp_2=q$
and $\delta \sim \delta^{-1/(p_1-1)}\eta(\Delta)^q$.
Consequently,
\[
\delta \sim \eta(\Delta)^{q-p}.
\]
Therefore,
\begin{equation*}
\begin{split}
\mathbb{E}[|X_{t}^{i}-X_{t}^{i,\Delta}|^{p}]
&\lesssim \int_{0}^{t}\mathbb{E}[|X_{s}^{i}-X_{s}^{i,\Delta}|^{p}]ds+\Delta^{p/2}h^p(\Delta)
+\eta(\Delta)^{q-p}.
\end{split}
\end{equation*}
Hence $\forall t \in [0,T]$,
\begin{equation*}
\sup_{0\leq t\leq T}\mathbb{E}[| X_{t}^{i}-X_{t}^{i,\Delta}|^{p}] \le \int_{0}^{t}\sup_{0\leq t\leq T}\mathbb{E}[|X_{s}^{i}-X_{s}^{i,\Delta}|^{p}]ds+C_{T}\Delta^{p/2}h^p(\Delta)
+C_{T}\eta(\Delta)^{q-p}.
\end{equation*}
Applying the Gr\"onwall inequality yields
\begin{equation*}
\begin{split}
\sup_{0\leq t\leq T}\mathbb{E}[|X_{t}^{i}-X_{t}^{i,\Delta}|^{p}]
\le C_{T}\Delta^{p/2}h^p(\Delta)+C_{T}\eta(\Delta)^{q-p}.
\end{split}
\end{equation*}

%Hence, one has
%\[
%\mathbb{E}[| X_{t}^{i}-X_{t}^{i,\Delta}|^{p}]
%\lesssim \Delta^{p/2}h(\Delta)^p+\eta(\Delta)^{q-p}.
%\]

Combining with the propagation of chaos results (Proposition \ref{pc}),  one therefore obtains the eventual estimate as listed in the statement.
\end{proof}

Now, we perform some discussions.
First of all, we note that if $b$ is globally Lipschitz and $a$ is one-sided Lipschitz, then $p_0$ can be arbitrarily large.
Then, $\eta(\Delta)$ can be chosen as polynomial of $\Delta$. For large enough $q$, the last term can be dropped. Hence, one has
\begin{crllr}
If $b$ is globally Lipschitz and $a$ is one-sided Lipschitz, then
\begin{equation}
 \sup_{t\le T}\mathbb{E}[|\bar{X}_{t}^{i}-X_{t}^{i,\Delta}|^{p}] \leq
 C\left\{\begin{array}{ll}N^{-1 / 2}+\Delta^{p/2}h(\Delta)^{p}, & \text { if } p>d/2, \\
  N^{-1 / 2} \log(1+N)+\Delta^{p/2}h(\Delta)^{p}, & \text { if } p=d/2, \\
 N^{-p/ d}+\Delta^{p/2} h(\Delta)^{p}, & \text { if }p\in[2,d/2).
 \end{array}\right.
\end{equation}
\end{crllr}
In this case, one may choose $h(\Delta)$ that grows slowly as $\Delta\to 0$ like $h(\Delta)\sim \Delta^{-\epsilon}$ for very small $\epsilon$.
This indicates that the order of the strong error can be arbitrarily close to one half.

However, if $p_0$ has an upper bound, choosing $h(\Delta)\sim \Delta^{-\epsilon}$ with small $\epsilon$ leads to big bound for the last term $\eta(\Delta)^{q/p-1}$.
In this case, one needs to choose suitable $h$ to optimize the rate.

\section{The random batch approximation for the interacting particle systems and its error analysis}\label{part4}

The computational cost for the simulation of the particle systems in \eqref{IN} and \eqref{INn} is high due to the $\cO(N^2)$ complexity per time step. In this section, we extend the recently proposed random batch method (RBM) \cite{jin2020,jin2022random,jin2021convergence,jin2022} to our particle systems and establish the approximation error estimate.

The Random Batch Method (RBM) is based on the so-called ``random mini-batch" idea, which is famous for its application in the stochastic gradient descent (SGD) \cite{robbins1951} algorithm for optimization in machine learning.  The key methodology of ``mini-batch'' is to find a cheap and unbiased random estimator using small subset of data/particles for the original quantity. How to design the random batch estimator clearly depend on the applications. In the RBM for particle systems, a random grouping strategy was proposed in \cite{jin2022random,jin2020,jin2021convergence}, while an importance sampling in the Fourier space was proposed for the Random Batch Ewald method for molecular dynamics in \cite{jin2022random,jin2021random}.  Note that for each single step, the random estimator has $\cO(1)$ error for the quantity. The key reason for the method to converge is the error cancelation in time. It is this type of  Law of Large Numbers (in time) that ensures convergence.  A difference of the RBM method from SGD is that the method is designed to dynamical properties of the systems, not just for equilibrium distribution.
%\tcr{For the geometric ergodicity and the long time behavior of the RBM  for interacting particle systems, see\cite{jin2022ergodicity}.   A fast potential energy splitting Markov chain Monte Carlo method was proposed in \cite{li2020}, which samples from the equilibrium distribution (Gibbs measure) of the particle system with singular interacting nuclei, and each step takes $\cO(1)$ time.}

In the original work \cite{jin2020},
\[
a(x, \mu)=f(x)+\int_{\R^d} k(x, y)\mu(dy), \quad b(x, \mu)\equiv \text{const}.
\]
Under this setting, the interacting particle systems becomes
\begin{gather}
d X_t^i=f(X_t^i)\, dt+\frac{1}{N-1}\sum_{j\neq i} k(X_t^i, X_t^j)dt+b\, dW_t^i.
\end{gather}
The initial values $X_0^i$ are i.i.d. sampled from $\mu_0$. Here, we have used $\frac{1}{N-1}\sum_{j\neq i}\delta_{X_t^i}$ to approximate $\mu$, and there is no significance difference regarding the propagation of chaos.

Suppose we aim to do simulation until time $T>0$. RBM does the following. Choose a batch size $P\ll N, P\ge 2$ that divides $N$. For the time step $\Delta$ and $t_m:=m\Delta$, on each time subinterval $[t_{m}, t_{m+1})$, there are two steps: (1) at time grid $t_{m}$, we divide the $N$ particles into $n:=N/P$ groups (batches) randomly; (2) the particles evolve with interaction inside the batches only.
\begin{lgrthm}[RBM]\label{rbm}
\begin{algorithmic}[1]
\noindent\FOR{$m \text{ in } 1: [T/\Delta]$}
\STATE Divide $\{1, 2, \ldots, N\}$ into $n=N/P$ batches randomly.
     \FOR{each batch  $\mathcal{C}_q$}
     \STATE Update $\tilde{X}^i_t$'s ($i\in \mathcal{C}_q$) by solving the following modified system with $t\in [t_{m}, t_{m+1})$.
     \begin{gather}\label{eq:rbmSDE}
            d \tilde{X}^i_t=f(\tilde{X}^i_t)\, dt+\frac{1}{P-1}\sum_{j\in \mathcal{C}_q,j\neq i} k(\tilde{X}_t^i, \tilde{X}_t^j)dt+b\, dW_t^i.
      \end{gather}
      \ENDFOR
\ENDFOR
\end{algorithmic}
\end{lgrthm}
In the case above, the dependence in $\mu$ is linear in the drift.
As proved in \cite{jin2020}, for a fixed configuration, the random batch approximation in the drift is unbiased and there is no restriction on batch size $P$.

Below we extend the random batch approximation to general coefficients where the diffusion coefficient could also depend on $\mu$ and they could be nonlinear in $\mu$. In general, if there is dependence of $\mu$ in $b$, the quadratic variation in the It\^o's formula would cause a contributing term of the order $1/P$ to the mean square error (or $\sqrt{1/P}$ to the error). Moreover, if the dependence on $\mu$ is nonlinear, the random approximation will no longer be unbiased for a fixed configuration. Consequently, it would also contribute an error that does not vanish as $\Delta\to 0$. All these would put restrictions on the batch size $P$. For example, $P\sim \Delta^{-1/2}$ to achieve half order for the mean square error. However, since the selection of batch size does not grow as the number of particles $N\to\infty$, the random batch approximation could still be beneficial.

We first consider a special case where the dependence on $\mu$ is linear to investigate how the random batch approximation in diffusion coefficient contributes to the error:
\begin{gather}\label{eq:diffusiondepend}
a(x,\mu)=f(x)+\int_{\mathbb{R}^{d}}k(x,y)\mu(dy),\quad
b(x,\mu)=\int_{\mathbb{R}^{d}}\sigma(x,y)\mu(dy).
\end{gather}
In particular, we investigate the bias introduced by the distribution dependence in the diffusion coefficient, which as we shall see,  is already very different from \cite{jin2020} regarding the application of RBM.

Secondly, we would consider drifts where the dependence on $\mu$ is nonlinear to investigate the error introduced by the random batch approximation:
\begin{gather}\label{eq:interactingnonlinear}
 a(x,\mu)=f(x)+A\left(\int k(x, y)\mu(dy)\right),\quad  b(x, \mu)\equiv \sigma(x).
\end{gather}
Here, $A: \R^r\to \R^d$ for some $r\ge 1$ is a nonlinear field.

Lastly, we perform discussion on general cases. In principle, the analysis would be similar to the cases discussed.  In this section, all the discussion on the random batch system will be performed by keeping the time continuous. Combining with the truncated scheme in the previous section, one will obtain the eventual method for numerical simulation.

\subsection{The particle system and the random batch approximation}

The interacting particle system \eqref{IN} for \eqref{eq:diffusiondepend} is given by
\begin{equation*}
dX^{i}_{t}=f(X^{i}_{t})\,dt+\frac{1}{N} \sum_{i=1}^{N} k(X^{i}_{t},X^{j}_{t})dt+ \frac{1}{N} \sum_{i=1}^{N} \sigma(X^{i}_{t},X^{j}_{t})dW^{i}_{t},
\end{equation*}
or
\begin{equation}\label{eq:interactinglinear}
dX^{i}_{t}=f(X^{i}_{t})\,dt+\frac{1}{N-1} \sum_{j\neq i}k(X^{i}_{t},X^{j}_{t})\,dt+ \frac{1}{N-1} \sum_{j\neq i}\sigma(X^{i}_{t},X^{j}_{t})dW^{i}_{t}.
\end{equation}
The initial values $X^{i}_{0} \in \R^d$ are the same as \eqref{IN}. The difference is whether we use $\frac{1}{N}\sum_j \delta_{X_t^j}$
or $\frac{1}{N-1}\sum_{j\neq i}\delta_{X_t^j}$ for the approximation.
Note that there is no big difference regarding the propagation of chaos while $1/(N-1)$ is more convenient in notation for the random batch system. We will focus on \eqref{eq:interactinglinear} then.
Similarly, for the case \eqref{eq:interactingnonlinear}, the interacting particle system is then given by
\begin{equation}\label{eq:particlenonlinear}
dX^{i}_{t}=f(X^{i}_{t})\,dt+A\left(\frac{1}{N-1} \sum_{j\neq i}k(X^{i}_{t},X^{j}_{t})\right)\,dt+ \sigma(X^i_t) dW^{i}_{t}.
\end{equation}

Recall that for the random batch approximation,  we divide the $N = nP$ particles into $n$ small batches with size $P\geq2$ randomly. Let the random batches be $\mathcal{C}_{q(i)}$, $q = 1, \ldots, n$.
We will use
\[
\cC:=\{\cC_1,\cdots, \cC_n \}
\]
to represent the division of random batches. Let $q(i)$ be the index $q$ such that $i\in \mathcal{C}_{q}$.
For the random batch approximation of \eqref{eq:interactinglinear}, we generate two random batches $\cC$ and $\cB$
at $t_{m}$ and run the following equation
\begin{equation}\label{eq:linearrbm}
d\tilde{X}^{i}_{t}=f(\tilde{X}_t^i)\,dt+\frac{1}{P-1} \sum_{j\in \cC_{q_1(i)},j\neq i} k(\tilde{X}^{i}_{t},\tilde{X}^{j}_{t})dt+ \frac{1}{P-1} \sum_{j\in \cB_{q_2(i)},j\neq i} \sigma(\tilde{X}^{i}_{t},\tilde{X}^{j}_{t})dW^{i}_{t}, \quad t\in [t_{m}, t_{m+1}),
\end{equation}
where $q_1\left(i\right)$ and  $q_2\left(i\right)$  are used to distinguish the batchmates in the two different batches  $\mathcal{B}$ and $\mathcal{C}$. And the initial values are $\tilde{X}^{i}_{0} =X_0^i$, the same as the particle system \eqref{IN}.

Correspondingly, the random batch approximation for \eqref{eq:particlenonlinear}
\begin{equation}\label{eq:nonlinearrbm}
d\tilde{X}^{i}_{t}=f(\tilde{X}^{i}_{t})\,dt+A\left(\frac{1}{P-1} \sum_{j\in \mathcal{C}_{q_{m}(i)},j\neq i}k(\tilde{X}^{i}_{t},\tilde{X}^{j}_{t})\right)\,dt+ \sigma(\tilde{X}^i_t) dW^{i}_{t},
\end{equation}
The initial values are $\tilde{X}^{i}_{0} =X_0^i$, the same as the particle system \eqref{IN}.

 Our experience tells us that if $N$ is not so small, there would be no significant difference if $\cC$ and $\cB$ are independent or the same. Hence,  for the convenience of the notations and analysis, we will always assume from here on that
\[
\cC=\cB.
\]
Namely, the random batches used for the drift and the diffusion are the same. To distinguish the batches at different time intervals, we denote the random division of batches at $t_{m}$ by
\[
\cC^{(m)}:=\{\cC_1^{(m)}, \cdots,
\cC_n^{(m)}\}.
\]

 Define the filtration
$\{\mathcal{F}_{m}\}_{m\geq 0}$ by
\begin{equation*}
\cF_{m}=\sigma\left( \cC^{(j)} : j \leq m \right).
\end{equation*}
This contains the information for the batches up to $t_m$. Also,
\[
\cF:=\sigma(\cup_{m}\cF_m)
\]
is the $\sigma$-algebra for how batches are generated at all time points.

\subsection{The linear dependence case}

In this section, we consider the special cases where the dependence on the distribution is linear and analyze the random batch approximation \eqref{eq:linearrbm}. For this purpose, we need stronger assumptions imposed on the coefficients than the one in Assumption \ref{assp3.2}. These conditions are still non-globally Lipschitz but the essential part is Lipschitz. Specifically,  we will investigate the problem with the following assumptions.
\begin{assumption}\label{ass:lipschitzass}
$f: \R^d\to \R^d$ is one-sided Lipschitz in the sense that
\begin{equation}\label{eq:lipschitzass}
(x-y)\cdot (f(x)-f(y))\le L|x-y|^2, \quad \forall x, y,
\end{equation}
and the derivatives have polynomial growth.
 $k:\mathbb{R}^{d}\times\mathbb{R}^{d}\to\mathbb{R}^{d}$ and $\sigma:\mathbb{R}^{d}\times\mathbb{R}^{d}\to \mathbb{R}^{d\times m'}$ are all Lipschitz continuous.
 \end{assumption}
\begin{rmrk}
The coefficient $b(x, \mu)=\int \sigma(x, y)\mu(dy)$ could also have a single term depending on $x$, like
$\sigma(x, y)=g(x)+\sigma_1(x, y)$.
In principle, $g(x)$ does not have to be globally Lipschitz. We assume the Lipschitz continuity of $\sigma$ simply for the convenience of discussion.
\end{rmrk}

\begin{prpstn}\label{prpstn:boundedness2}
Let Assumption \ref{ass:lipschitzass} hold and $\mu_0\in \cP_q$ with $q\ge 2$. Let $ X^j_t$ be the solution to \eqref{eq:interactinglinear} and $\tilde{X}^j_t$ be the solution to \eqref{eq:linearrbm} for $j=1,\cdots,N$. Then conditioning on the sequence of random batches, one has
\begin{equation}
\sup_{0\le t\le T}\max_{1\le j\le N}\E[|\tilde{X}^j_t|^q  | \cF ]
\le C_q.
\end{equation}
Similarly,
\begin{equation}
\sup _{0\leq t\leq T}\left(\mathbb{E}[\left|X^{i}_{t}\right|^{q}]+\mathbb{E}[|\tilde{X}^{i}_{t}|^{q}]\right)=\sup _{0\leq t\leq T}\left(\mathbb{E}[\left|X^{1}_{t}\right|^{q}]+\mathbb{E}[|\tilde{X}^{1}_{t}|^{q}]\right) \leq C_{q,T}, \forall i.
\end{equation}
In both estimates, the constant $C_q(T)$ is independent of $N$ and the sequence of batches.
\end{prpstn}

\begin{proof}
We only show the moment control for $\tilde{X}^i_t$ because the moment control for $X^i_t$ is similar. The proof is essentially
the same as \cite[Lemma 3.4]{jin2022}.  Our approach is to show the boundedness of the moments for any given sequences of the batches, to avoid the interplay between randomness of the batches and the process $\tilde{X}_t$.

For $t\in [t_{m}, t_{m+1})$, by It\^o's formula one has
\begin{multline*}
\frac{d}{dt}\E[|\tilde{X}^i_t|^q | \cF ]
=q\E\left[|\tilde{X}^i_t|^{q-2}\tilde{X}^i_t
\cdot\Big(f(\tilde{X}^i_t)+\frac{1}{P-1}\sum_{j\in \mathcal{C}_{q_{m}(i)},j\neq i}k(\tilde{X}^i_t, \tilde{X}^j_t)\Big) \Bigg| \cF \right]+\\
\frac{1}{2}q \E\left[|\tilde{X}^i_t|^{q-2}\Big(I+(q-2)\frac{\tilde{X}^i_t\otimes \tilde{X}^i_t}
{|\tilde{X}^i_t|^2} \Big): \Big(\frac{1}{P-1}\sum_{j\in \mathcal{C}_{q_{m}(i)},j\neq i}\sigma(\tilde{X}^i_t,\tilde{X}^j_t)\Big) \Big(\frac{1}{P-1}\sum_{j\in \mathcal{C}_{q_{m}(i)},j\neq i}\sigma(\tilde{X}^i_t,\tilde{X}^j_t)\Big)^T \Bigg| \cF \right]\\
=: M_1+M_2.
\end{multline*}
By Assumption \ref{ass:lipschitzass},
\[
\begin{split}
& \tilde{X}^i_t\cdot f(\tilde{X}^i_t)\le L|\tilde{X}^i_t|^2+|f(0)||\tilde{X}^i_t|,
\\
&|\tilde{X}^i_t|^{q-2}\tilde{X}^i_t\cdot k(\tilde{X}^i_t, \tilde{X}^j_t))
\le L((2-q^{-1})|\tilde{X}^i_t|^q+q^{-1}|\tilde{X}^j_t|^q)
+C|\tilde{X}^i_t|^{q-1}.
\end{split}
\]
One thus has
\[
M_1\lesssim 1+\E[|\tilde{X}^i_t|^q | \cF ]+\frac{1}{P-1}
\sum_{j\in\mathcal{C}_{q_{m}(i)},j\neq i}\E(|\tilde{X}^j_t|^q  | \cF )
\lesssim 1+\max_{1\le j\le N}\E(|\tilde{X}^j_t|^q  | \cF ).
\]
Similarly,
\[
M_2\le
q(q-1)\frac{1}{(P-1)^2}\sum_{j\in \mathcal{C}_{q_{m}(i)}, j'\in\mathcal{C}_{q_{m}(i)},j\neq i, j'\neq i}
\E( |\tilde{X}^i_t|^{q-2}\|\sigma(\tilde{X}^i_t,\tilde{X}^j_t)\sigma(\tilde{X}^i_t,\tilde{X}^{j'}_t)^T\|  | \cF ).
\]
Again, using the Lipschitz continuity of $\sigma$ and the Young's inequality like
\[
|\tilde{X}^i_t|^{q-2} |\tilde{X}^j_t||\tilde{X}^{j'}_t|
\le \frac{q-2}{q}|\tilde{X}^i_t|^{q}+\frac{1}{q}|\tilde{X}^j_t|^q
+\frac{1}{q}|\tilde{X}^{j'}_t|^q,
\]
one will get the same bound as in $M_1$. Consequently, by gluing together all the estimates on different subintervals of the form
$[t_{m}, t_{m+1})$,  one has
\[
\E[|\tilde{X}^i_t|^q | \cF]
\le \E[|\tilde{X}^i_0|^q | \cF ]
+C\int_0^t (1+\max_{1\le j\le N}\E[|\tilde{X}^j_s|^q  | \cF ])\,ds,
\]
for any $i$. Hence, defining
\[
a(t):=\max_{1\le j\le N}\E[|\tilde{X}^j_t|^q  | \cF ],
\]
one has the inequality
\[
a(t)\le a(0)+C\int_0^t a(s)\,ds.
\]
Applying Gr\"onwall's inequality yields the result for $\E[|\tilde{X}^j_t|^q  | \cF ]$. Then, taking expectation over the random batches, the result for the moments also follows.
\end{proof}

Below, for a given configuration, $\mathtt{x}:=(x^1, x^2, \cdots, x^N)
\in \R^{Nd}$. Define the random deviation in a random batch approximation, with interaction kernel $k$:
\begin{gather}
\chi_i(\mathtt{x}; k):= \frac{1}{P-1} \sum_{j\in \mathcal{C}_{q_{m}(i)},j\neq i} k(x^i, x^j)-\frac{1}{N-1} \sum_{j\neq i}k(x^i,x^j).
\end{gather}

The following lemma, shown in  \cite[Lemma 3.1]{jin2020}, gives the basic consistency and stability features of the random batch method.
\begin{lmm}\label{lmm:consistencyofrbm}
The random deviation satisfies
\begin{equation*}
\E [\chi_i(\mathtt{x}; k)]=0
\end{equation*}
and
\begin{gather}
\var \chi_i(\mathtt{x}; k)=\E[\|\chi_i(\mathtt{x}; k)\|^2]=\left(\frac{1}{P-1}-\frac{1}{N-1}\right) \Lambda_i(\mathtt{x};k),
\end{gather}
where
\[
\Lambda_i(\mathtt{x};k) := \frac{1}{N-2} \sum_{j: j \neq i}\Big\|k(x^i, x^j)-\frac{1}{N-1} \sum_{j': j' \neq i} k(x^i, x^{j'}) \Big\|^{2}.
\]
\end{lmm}

We remark that if $k$ is a matrix-valued function, the norm in this lemma is the Hilbert-Schmidt norm.
Let us briefly understand how the random batch approximation in the diffusion coefficient affects the error.
If we define
\begin{gather}\label{eq:errorprocess}
Z^i_t := \tilde{X}^i_t-X^i_t.
\end{gather}
Then, one finds that
\begin{equation*}
\begin{split}
dZ^i_t&= (f(\tilde{X}^i_t)-f(X^i_t))\,dt
+\left(\frac{1}{P-1}\sum_{j\in \mathcal{C}_{q_{m}(i)},j\neq i}k(\tilde{X}^i_t, \tilde{X}^j_t)
-\frac{1}{N-1}\sum_{j\neq i}k(X^i_t, X^j_t)   \right)\,dt\\
&\quad+\left(\frac{1}{P-1}\sum_{j\in \mathcal{C}_{q_{m}(i)},j\neq i}\sigma(\tilde{X}^i_t, \tilde{X}^j_t)
-\frac{1}{N-1}\sum_{j\neq i}\sigma(X^i_t, X^j_t)   \right)dW^i_t.
\end{split}
\end{equation*}
Intuitively, if $X^i_t$ is close to $\tilde{X}^i_t$, then the leading error would be
\[
dZ^i_t \sim \chi_{i}(\bfX;k)dt+\chi_{i}(\bfX; \sigma)dW^i_t,
\]
where $\mathbf{X}:=(X^1, X^2,\cdots, X^N)\in \mathbb{R}^{Nd}$. %where the batches in $\chi_1$ and $\chi_2$ are assumed to be independent.
 Then, by It\^o's  formula,
\[
d\E [|Z^i_t|^2]\sim 2 \E [Z^i_t\cdot \chi_{i}(\bfX;k)]\,dt+\E [|\chi_{i}(\bfX;\sigma)|^2]\,dt.
\]
The batch at $t_{m}$ is independent of $Z^i_{t_{m}}$ so conditioning on $Z^i_{t_{m}}$, the expectation of $\chi_i$ is zero, then the term $Z^i_t\cdot \chi_i$ would contribute
a term that vanishes as $\Delta\to 0$. The second term, however, is of order $1/P$. This would not vanish as $\Delta\to 0$ if we fix $P$. This term  is the quadratic variation due to the Brownian motion.  Clearly, when there is the distribution dependence in the diffusion coefficient, the effect of random batch would be very different from the ones discussed in \cite{jin2020}.

The following result gives the basic estimate for the approximation error.
\begin{thrm}\label{thm:batch1}
Let Assumption \ref{ass:lipschitzass} hold. $\bfX(t):=(X^1_t, \cdots, X^N_t)$ are the solution of \eqref{eq:interactinglinear} and $\tilde{\bfX}(t):=(\tilde{X}^1_t, \cdots, \tilde{X}^N_t)$ are the solution of \eqref{eq:linearrbm}. Then, the error process defined in \eqref{eq:errorprocess} satisfies
\begin{equation}
\sup_{0\le t\le T}\max_i \E [|Z^i_t|^2]\le C\left(\Delta^2+\frac{ \Lambda(k) \Delta}{P}
+\frac{\Lambda(\sigma)}{P}\right),
\end{equation}
where $\Lambda(k):=\sup\limits_{0\le t\le T}\E [\Lambda_i(\bfX(t); k)]<\infty$
and $\Lambda(\sigma):=\sup\limits_{0\le t\le T} \E [\Lambda_i(\bfX(t); \sigma)]<\infty$ are independent of $i$ and can be bounded independent of $N$.
In particular,
\begin{gather}\label{eq:mainestimate}
\sup_{0\le t\le T} \E [|Z^i_t|^2]\le \begin{cases}
C_T(\frac{\Delta}{P}+\Delta^2), & \sigma(x, y)\equiv \sigma(x),\\
C_T(\frac{1}{P}+\Delta^2), & \text{otherwise}.
\end{cases}
\end{gather}
\end{thrm}

Clearly, when $\sigma(x, y)\equiv \sigma(x)$, the result \eqref{eq:mainestimate}  is the same
as in \cite{jin2020}, though they assumed boundedness of the kernel $k(x, y)$.  For general $\sigma$ that may depend on $y$, a fixed batch size $P$ clearly would not give convergence. Instead, we need to take $P$ large to get convergence.
\begin{crllr}\label{crll}
If one takes $P=\min(\Delta^{-\beta}, N)$ with $\beta\le 1$, which stays finite as $N\to\infty$, the error is controlled as
\begin{equation}\label{eq:beta}
\sqrt{\E [|Z^i_t|^2]}\le \begin{cases}
C\Delta^{(\beta+1)/2}, & \sigma(x, y)\equiv \sigma(x),\\
C\Delta^{\beta/2}, & \text{otherwise}.
\end{cases}
\end{equation}
\end{crllr}

To prove Theorem \ref{thm:batch1}, we need the following symmetry fact.
\begin{lmm}\label{lmm:symm}
Let $Z^i_t$ be given in \eqref{eq:errorprocess}, with $X_0^i$ being i.i.d. sampled. Let $\varphi(\cdot)$ be a test function that has at most polynomial growth. Then, for $[t_{m}, t_{m+1})$, one has
\[
\E\left[\E\Big[\frac{1}{P}\sum_{j\in \mathcal{C}_{q_{m}(i)}} \varphi(Z^j_t) | \cF_{m}\Big]\right]
=\frac{1}{N}\sum_{j'=1}^N\E[\varphi(Z^{j'}_t)]= \E[\varphi(Z^j_t)], \quad \forall j\in \{1,\cdots, N\}.
\]
Similarly, for any particle $i$,
\[
\E\left[\E\Big[\frac{1}{P-1}\sum_{j\in \mathcal{C}_{q_{m}(i)}, j\neq i} \varphi(Z^j_t) | \cF_{m}\Big]\right]
=\frac{1}{N-1}\sum_{j': j'\neq i} \E[\varphi(Z^{j'}_t)]=\E[\varphi(Z^j_t)], \quad \forall j\in \{1,\cdots, N\}.
\]
\end{lmm}

\begin{proof}
 Note that $\E[\E[\cdot | \cF_m]]=\E[\E[\cdot | \sigma(\cC^{(m)})]]$. Then, by definition,
\begin{equation*}
\begin{split}
&\E\left[\E\left[\frac{1}{P}\sum\limits_{i\in \mathcal{C}_{q_{m}(i)}} \varphi(Z^i_t) | \sigma(\cC^{(m)}) \right]\right]\\
&=\sum_{C_1,\cdots, C_n}
\P(\cC_1=C_1,\cdots, \cC_n=C_n)\frac{1}{P}\sum_{j'=1}^{N}
I_{j'\in \mathcal{C}_{q_{m}(i)}}
\E[ \varphi(Z^{j'}_t) | \cC_1=C_1,\cdots, \cC_n=C_n]\\
&
=\frac{1}{P}\sum_{j'=1}^{N}
\sum_{C_1,\cdots, C_n}
\P(\cC_1=C_1,\cdots, \cC_n=C_n)I_{j'\in \mathcal{C}_{q_{m}(i)}}
\E[ \varphi(Z^{j'}_t) | \cC_1=C_1,\cdots, \cC_n=C_n].
\end{split}
\end{equation*}
Here, we make of the symmetry.
The point is that $\E[ \varphi(Z^{j'}(t)) | \cC_1=C_1,\cdots, \cC_n=C_n]$ does not depend on the batch (the values $C_1,\cdots, C_n$). In other words, no matter what its batchmates are, this conditional expectation would be the same because
 the joint distribution of the particles is symmetric at $t_{m}$.
 Hence,
 \begin{equation*}
 \begin{split}
&\sum_{C_1,\cdots, C_n}
\P(\cC_1=C_1,\cdots, \cC_n=C_n)I_{j'\in C_q}
\E[ \varphi(Z^{j'}_t) | \cC_1=C_1,\cdots, \cC_n=C_n]\\
&=\P(j'\in \mathcal{C}_{q_{m}(i)}) \E[\varphi(Z^{j'}_t)]=\frac{P}{N}\E[\varphi(Z^{j'}_t)].
\end{split}
 \end{equation*}
The second claim is nearly the same. The only difference is that the batch $\mathcal{C}_{q_{m}(i)}$ used will vary. Let $\mathbf{1}_{ij'}$ be the indicator that particles $i$ and $j'$ are in the same batch. Then,
\begin{equation*}
\begin{split}
&\sum_{C_1,\cdots, C_n}
\P(\cC_1=C_1,\cdots, \cC_n=C_n)\frac{1}{P-1}\sum_{j'=1}^{N}
\mathbf{1}_{ij'}
\E[ \varphi(Z^{j'}_t) | \cC_1=C_1,\cdots, \cC_n=C_n]\\
&
=\frac{1}{P-1}\sum_{j'=1}^{N}\P(\mathbf{1}_{ij'}=1)
\sum_{C_1,\cdots, C_n}
\P(\cC_1=C_1,\cdots, \cC_n=C_n)
\E[ \varphi(Z^{j'}_t) | \cC_1=C_1,\cdots, \cC_n=C_n].
\end{split}
\end{equation*}
Applying \cite[Lemma 3.2]{jin2021convergence} gives the result.
\end{proof}

\begin{rmrk}
This lemma is different from \cite[Lemma 3.4]{jin2021convergence}, where the summands are required to be independent from the random batches. Here, we do not require the independence, but use the symmetry. Nevertheless the goal is the same: to deal with the random sum.
\end{rmrk}

We rearrange the terms in $dZ^i_t$ as
\begin{align}
dZ^i_t&=\left(\int_0^{1}\nabla f\left(X^i_t+\tau\left(\tilde{X}^i_t-X^i_t\right)\right)d\tau\right)\cdot Z^i_t
\notag\\&\quad+\frac{1}{P-1}\sum_{j\in \mathcal{C}_{q_{m}(i)},j\neq i}\left(k\left(\tilde{X}^i_t, \tilde{X}^j_t\right) - k\left(X^i_t, X^j_t\right)\right)dt
+\chi_i\left(\bfX(t);k\right)dt\notag\\&
\quad+\frac{1}{P-1}\sum_{j\in \mathcal{C}_{q_{m}(i)},j\neq i}\left(\sigma\left(\tilde{X}^i_t, \tilde{X}^j_t\right) - \sigma\left(X^i_t, X^j_t\right)\right)dW^i_t
+\chi_i\left(\bfX(t); \sigma\right)dW^i_t.
\notag\end{align}
Now, we can prove the main result.

\begin{proof}[Proof of Theorem \ref{thm:batch1}]
For $t\in [t_m, t_{m+1})$, by It\^o's formula, one has
\begin{align}
 \frac{d}{dt}&\E\left[ \left|Z^i_t\right|^2 | \cF_{m} \right]
\notag\\&=\Bigg(2 \E\bigg[ Z^i_t\cdot \left(f\left(X^i_t\right)-f(\tilde{X}^i_t)\right) \bigg| \cF_{m} \bigg]
+2\E\bigg[\frac{1}{P-1}\sum_{j\in \mathcal{C}_{q_{m}(i)},j\neq i}
Z^i_t \cdot \left(k\left(\tilde{X}^i_t, \tilde{X}^j_t\right)-k\left(X^i_t, X^j_t\right)\right) \bigg| \cF_{m} \bigg]\Bigg)
\notag\\
& \quad+2\E\left[ Z^i_t\cdot \chi_{i}\left(\bfX(t); k\right) |  \cF_{m} \right]  +\E\left[\Sigma\Sigma^{\scriptscriptstyle T} | \cF_{m} \right]
\notag\\&=S_1+S_2+S_3,\notag
\end{align}
where
%\begin{align}
%S_1&=2 \E\bigg( Z^i_t\cdot \left(f\left(X^i_t\right)-f(\tilde{X}^i_t)\right) \bigg| \cF_{m} \bigg)
%\notag\\&\quad+2\E\bigg(\frac{1}{P-1}\sum_{j\in \mathcal{C}_{q_{m}(i)},j\neq i}
%Z^i_t \cdot \left(k\left(\tilde{X}^i_t, \tilde{X}^j_t\right)-k\left(X^i_t, X^j_t\right)\right) \bigg| \cF_{m} \bigg),
%\notag\\S_2&=2\E\left( Z^i_t\cdot \chi_{i}\left(\bfX(t); k\right) |  \cF_{m} \right),
%\notag\\S_3&=\E\left(\Sigma\Sigma^{\scriptscriptstyle T} | \cF_{m} \right).\notag
%\end{align}
%and
\[
\Sigma=\frac{1}{P-1}\sum_{j\in \mathcal{C}_{q_{m}(i)},j\neq i}\sigma(\tilde{X}^i_t, \tilde{X}^j_t)
-\frac{1}{N-1}\sum_{j\neq i}\sigma(X^i_t, X^j_t).
\]

For $S_1$ term, using the one-sided Lipschitz condition of $f$ and the Lipschitz continuity of $k$, one has
\begin{equation*}
\begin{split}
S_1  & \lesssim  \E\left[ \left|Z^i_t\right|^2 | \cF_{m}  \right]
+\E\bigg[\frac{1}{P-1}\sum_{j\in \mathcal{C}_{q_{m}(i)},j\neq i}|Z^i_t|\left(|Z^i_t|+|Z^j_t|\right)\Big| \cF_{m} \bigg] \\
& \lesssim  \E\left[ \left|Z^i_t\right|^2 | \cF_{m}\right] +\E\bigg[\frac{1}{P-1}\sum_{j\in \mathcal{C}_{q_{m}(i)}, j\neq i}\left|Z^j_t\right|^2 \Big| \cF_{m}\bigg].
\end{split}
\end{equation*}

The $S_3$ term can be estimated similarly:
\begin{gather*}
\begin{split}
S_3  &\le  \E\left[\left\|\Sigma\right\|^2 | \cF_{m}\right]
\\&\le 2\E\bigg[\Big\|  \frac{1}{P-1}\sum_{j\in \mathcal{C}_{q_{m}(i)},j\neq i}\sigma(\tilde{X}^i_t, \tilde{X}^j_t)  -\frac{1}{P-1}\sum_{j\in \mathcal{C}_{q_{m}(i)},j\neq i}\sigma\left(X^i_t, X^j_t\right) \Big\|^2 \Big|\cF_{m} \bigg]
+ 2\E\left[\|  \chi_{i}\left(\bfX(t); \sigma\right) \|^2 | \cF_{m} \right]
\\&=: S_{31}+S_{32}.
\end{split}
\end{gather*}
Here, $S_{32}$ is the main local error term arising from $S_3$, which will not vanish as $\Delta\to 0$. The $S_{31}$ term can be estimated using the Lipschitz continuity of $\sigma$ so that
\begin{equation*}
\begin{split}
S_{31} & \lesssim \E\left[\left|Z^i_t\right|^2 | \cF_{m}\right]
+\E\bigg[\Big|\frac{1}{P-1}\sum_{j\in \mathcal{C}_{q_{m}(i)},j\neq i}\left|Z^j_t\right|  \Big|^2 \Big|\cF_{m}\bigg] \\
&\le  \E\left[\left|Z^i_t\right|^2 | \cF_{m}\right]+\E\bigg[\frac{1}{(P-1)^2}\sum_{j\in \mathcal{C}_{q_{m}(i)}, j'\in \mathcal{C}_{q_{m}(i)},j,j'\neq i}\frac{1}{2}\left(\left|Z^j_t\right|^2 +\left|Z^{j'}_t\right|^2 \right)\Big|\cF_{m}\bigg]\\
&=\E\left[\left|Z^i_t\right|^2 | \cF_{m}\right]+\E\bigg[\frac{1}{P-1}\sum_{j\in \mathcal{C}_{q_{m}(i)},j\neq i} \left|Z^j_t\right|^2 \Big|\cF_{m}\bigg].
\end{split}
\end{equation*}

Now, for $S_2$,  due to Lemma \ref{lmm:consistencyofrbm} and the Lipschitz continuity of $k$, we break this term into the following.
\allowdisplaybreaks[4]
\begin{align}
S_2&\lesssim \left\{\E\left[Z^i_{t_{m}}\cdot \chi_{i}\left(\bfX(t); k\right)\Big|\cF_{m} \right]\right\}
 +\left\{\E\bigg[\left|\chi_{i}\left(\bfX(t); k\right)\right|\int_{t_{m}}^t\left(1+\left|\tilde{X}^i_s\right|^{\gamma}+\left|X^i_s\right|^{\gamma}\right)\left|Z^i_s\right|ds \Big|\cF_{m}\bigg]\right.\notag\\
&\quad\left.+\E\bigg[\left|\chi_{i}\left(\bfX(t); k\right)\right|  \int_{t_{m}}^t\bigg(\left|Z^i_s\right|
+\frac{1}{P-1}\sum_{j\in\mathcal{C}_{q_{m}(i)},j\neq i}|Z^j_s|\bigg)ds \Big|\cF_{m}\bigg]\right\}\notag\\
&\quad+ \left\{\E\bigg[\chi_i\left(\bfX(t); k\right)\int_{t_{m}}^t \frac{1}{P-1}\sum_{j\in \mathcal{C}_{q_{m}(i)},j\neq i}
\left(\sigma\left(\tilde{X}^i_s, \tilde{X}^j_s\right)-\sigma\left(X^i_s, X^j_s\right)\right)dW^i_s  \Big| \cF_{m}\bigg]\right\}\notag\\
&\quad+ \left\{\int_{t_{m}}^t\E\Big[ \chi_{i}\left(\bfX(s);k\right)
\cdot\chi_i\left(\bfX(t); k\right) | \cF_{m}\Big]ds
+\E\bigg[ \chi_{i}\left(\bfX(t); k\right)\int_{t_{m}}^t \chi_i\left(\bfX(s); \sigma\right) dW^i_{s}    \Big|\cF_{m}\bigg] \right\} \notag\\
 &=:S_{20}+S_{21}+S_{22}+S_{23}.\notag
\end{align}
The terms around the brackets on the right side are called $S_{2,i}$, $i=0, 1, 2, 3$ respectively for the convenience of notations.
%where
%\begin{align}
%S_{20}&=\E\left(Z^i_{t_{m}}\cdot \chi_{i}\left(\bfX(t); k\right)|\cF_{m} \right);\notag
%\\S_{21}&=\E\bigg(\Big(\chi_{i}\left(\bfX(t); k\right)\Big)\int_{t_{m}}^t\left(1+\left|\tilde{X}^i_s\right|^{\gamma}+\left|X^i_s\right|^{\gamma}\right)\left|Z^i_s\right|ds |\cF_{m}\bigg)\notag\\
%&\quad+\E\bigg(\left|\chi_{i}\left(\bfX(t); k\right)\right|  \int_{t_{m}}^t\bigg(\left|Z^i_s\right|
%+\frac{1}{P-1}\sum_{j\in\mathcal{C}_{q_{m}(i)},j\neq i}|Z^j_s|\bigg)ds \bigg|\cF_{m}\bigg);\notag
%\\S_{22}&= \E\bigg( \chi_i\left(\bfX(t); k\right)\int_{t_{m}}^t \frac{1}{P-1}\sum_{j\in \mathcal{C}_{q_{m}(i)},j\neq i}
%\left(\sigma\left(\tilde{X}^i_s, \tilde{X}^j_s\right)-\sigma\left(X^i_s, X^j_s\right)\right)dW^i_s  \Big| \cF_{m}\bigg);\notag
%\\S_{23}&=\int_{t_{m}}^t\E\Big( \chi_{i}\left(\bfX(s);k\right)
%\cdot\chi_i\left(\bfX(t); k\right) | \cF_{m}\Big)ds
%\notag\\&\quad+\E\left( \chi_{i}\left(\bfX(t); k\right)\int_{t_{m}}^t \chi_i\left(\bfX(s); \sigma\right) dW^i_{s}    |\cF_{m} \right).\notag
%\end{align}
We will leave $S_{20}$ as it is. Consider $S_{21}$ first.  The first term in $S_{21}$ can be controlled using Proposition \ref{prpstn:boundedness2} by
\begin{align}
 \int_{t_{m}}^t &\Big(\E\left(\chi_{i}\left(\bfX(t); k\right)^4 |\cF_{m}\right)\Big)^{1/4}
\bigg(\E\left(\left(1+\left|\tilde{X}^i_s\right|^{\gamma}+\left|X^i_s\right|^{\gamma}\right)^4 |\cF_{m}\right)\bigg)^{1/4}
\notag\\&\times\bigg(\E\left(\left|Z^i_s\right|^2 |\cF_{m}\right)\bigg)^{1/2}ds
\le C\int_{t_{m}}^t\bigg(\E\left(\left|Z^i_s\right|^2 |\cF_{m}\right)\bigg)^{1/2}ds.\notag
\end{align}
The second term in $S_{21}$ can be estimated similarly by
$$
 \int_{t_{m}}^t\bigg(\E\left[\left|Z^i_s\right|^2 |\cF_{m}\right]\bigg)^{1/2}ds
+ \int_{t_{m}}^t\bigg(\E\bigg[\frac{1}{P-1}\sum_{j\in \mathcal{C}_{q_{m}(i)},j\neq i}\left|Z^j_s\right|^2\Big|\cF_{m}\bigg]\bigg)^{1/2}ds.
$$
The $S_{22}$ term can be estimated using the Lipschitz continuity of $\sigma$, Proposition \ref{prpstn:boundedness2} and the It\^o's isometry by
\begin{equation*}
\begin{split}
\int_{t_{m}}^t &\E\Bigg[ \bigg(\left|Z^i_s\right|+\frac{1}{P-1}\sum_{j\in \mathcal{C}_{q_{m}(i)},j\neq i}
\left|Z^j_s\right|\bigg)^2 \Big|\cF_{m}\Bigg]^{1/2}ds \\&\le
\int_{t_{m}}^t\Big(\E\left[|Z^i_s|^2 |\cF_{m}\right]\Big)^{1/2}ds
+ \int_{t_{m}}^t\bigg(\E\bigg[\frac{1}{P-1}\sum_{j\in \mathcal{C}_{q_{m}(i)},j\neq i}\left|Z^j_s\right|^2\Big|\cF_{m}\bigg]\bigg)^{1/2}ds.
\end{split}
\end{equation*}
Here, $S_{23}$ term contains essentially the local error term for random batch approximation arising from $S_2$. We may expect that this tends to zero as $\Delta\to 0$.  This term can be simply estimated by
\begin{align}
S_{23}&\le \frac{1}{2}\E\left[\left|\chi_i\left(\bfX(t); k\right)\right|^2 |\cF_{m}\right]\Delta
+\frac{1}{2}\int_{t_{m}}^t \E\left[\left|\chi_i\left(\bfX(s); k\right)\right|^2 |\cF_{m}\right]ds
\notag\\&\quad
+\E\left[\left|\chi_i\left(\bfX(t); k\right)\right|^2 |\cF_{m}\right]^{1/2}\left(\int_{t_{m}}^t \E\left[\left\|\chi_i\left(\bfX(s); \sigma\right)\right\|^2 |\cF_{m}\right]ds \right)^{1/2}.\notag
\end{align}
Combining all these estimates together, we find
\begin{align}
 \E&\left[\left|Z^i_t\right|^2|\cF_{m}\right]
\notag\\& \le \E\left[\left|Z^{i}_{t_{m}}\right|^2|\cF_{m}\right]+C\int_{t_{m}}^t \E\left[\left|Z^i_s\right|^2|\cF_{m}\right]ds
+C\int_{t_{m}}^t\E\bigg[\frac{1}{P-1}\sum_{j\in \mathcal{C}_{q_{m}(i)},j\neq i} \left|Z^j_s\right|^2 \Big|\cF_{m}\bigg]ds
\notag\\&\quad+C\int_{t_{m}}^t \int_{t_{m}}^s\left(\E\left[\left|Z^i_\tau\right|^2 |\cF_{m}\right]\right)^{1/2}d\tau ds+\int_{t_{m}}^t \int_{t_{m}}^s\bigg(\E\bigg[\frac{1}{P-1}\sum_{j\in \mathcal{C}_{q_{m}(i)},j\neq i}\left|Z^j_\tau\right|^2\Big|\cF_{m}\bigg]\bigg)^{1/2}d\tau ds \notag\\
&\quad+\int_{t_{m}}^t \E\left[Z^i_{t_{m}}\cdot \chi_{i}\left(\bfX(s); k\right) |\cF_{m}\right]ds+ \frac{\Delta}{2}\int_{t_{m}}^t\E\left[\left|\chi_i\left(\bfX(s); k\right)\right|^2 |\cF_{m}\right]ds
+\frac{1}{2}\int_{t_{m}}^t\int_{t_{m}}^s \E\left[\left|\chi_i(\bfX(\tau); k)\right|^2 |\cF_{m}\right]d\tau ds
\notag\\
&\quad+\int_{t_{m}}^t \E\Big[\left|\chi_i\left(\bfX(s); k\right)\right|^2 |\cF_{m}\Big]^{1/2}\Big(\int_{t_{m}}^s \E\left[\left\|\chi_i\left(\bfX(\tau); \sigma\right)\right\|^2 |\cF_{m}\right]d\tau\Big)^{1/2} ds
+2\int_{t_{m}}^t\E\left[\left\|  \chi_{i}\left(\bfX(s); \sigma\right) \right\|^2 \Big| \cF_{m} \right]ds.\notag
\end{align}
The double integrals can be estimated, for example, by
\begin{align}
\int_{t_{m}}^t \int_{t_{m}}^s\left(\E\left[\left|Z^i_\tau\right|^2 |\cF_{m}\right]\right)^{1/2}d\tau ds
=\int_{t_{m}}^t (t-\tau)\left(\E\left[\left|Z^i_\tau\right|^2 |\cF_{m}\right]\right)^{1/2}d\tau
\le \frac{\Delta^3}{6}+\frac{1}{2}\int_{t_{m}}^t\E\left[\left|Z^i_\tau\right|^2 |\cF_{m}\right]d\tau,\notag
\end{align}
and
\begin{align}
& \int_{t_{m}}^t \E\left[\left|\chi_i\left(\bfX(s); k\right)\right|^2 |\cF_{m}\right]^{1/2}\Big(\int_{t_{m}}^s \E\left[\left\|\chi_i\left(\bfX(\tau); \sigma\right)\right\|^2 |\cF_{m}\right]d\tau\Big)^{1/2} ds \notag\\
&\le  \frac{\Delta}{2}\int_{t_{m}}^t\E\left[\left|\chi_i\left(\bfX(s); k\right)\right|^2 |\cF_{m}\right]ds
+\frac{1}{2\Delta}\int_{t_{m}}^t(t-\tau) \E\left[\left\|\chi_i\left(\bfX(\tau); \sigma\right)\right\|^2 |\cF_{m}\right]d\tau.\notag
\end{align}
One actually has
\begin{align}
  \E&\left[\left|Z^i_t\right|^2|\cF_{m}\right]
\notag\\&\le \E\left[\left|Z^{i}_{t_{m}}\right|^2|\cF_{m}\right]+C\int_{t_{m}}^t \E\left[\left|Z^i_s\right|^2|\cF_{m}\right]ds
+C\int_{t_{m}}^t\E\bigg[\frac{1}{P-1}\sum_{j\in \mathcal{C}_{q_{m}(i)},j\neq i} \left|Z^j_s\right|^2 \bigg|\cF_{m}\bigg]ds
\notag\\&\quad+\int_{t_{m}}^t\E\left[Z^i_{t_{m}}\cdot \chi_{i}\left(\bfX(s); k\right) |\cF_{m}\right]ds
+C\Delta^3+\frac{3\Delta}{2}\int_{t_{m}}^t\E\left[\left|\chi_i\left(\bfX(s); k\right)\right|^2 |\cF_{m}\right]ds
+\frac{5}{2}\int_{t_{m}}^t\E\left[\left\|  \chi_{i}\left(\bfX(s); \sigma\right) \right\|^2 \Big| \cF_{m}\right]ds.\notag
\end{align}

%Now, we take expection over $\cF_m$.
Due to Lemma \ref{lmm:consistencyofrbm}, $ \int_{t_{m}}^t\E(Z^i_{t_{m}}\cdot \chi_{i}(\bfX(s); k) |\cF_{m})ds=0$. Regarding the term
$$
\int_{t_{m}}^t\E\left[ \frac{1}{P-1}\sum_{j\in \mathcal{C}_{q_{m}(i)},j\neq i} \left|Z^j_s\right|^2\bigg|\cF_{m}\right]ds
$$
we note that the batch $\cC^{(m)}$ is not independent of $|Z^j_t|$. To treat this, \cite{jin2020,jin2021convergence} considered the difference $|Z^j_t|-|Z^j_{t_{m}}|$.  Here, we make use of the symmetry. By Lemma \ref{lmm:symm}, taking expectation on both sides, one has
\begin{equation*}
 \E\left[\left|Z^i_t\right|^2\right]  \le \E\left[\left|Z^{i}_{t_{m}}\right|^2\right]+C\int_{t_{m}}^t \E\left[\left|Z^i_{s}\right|^2\right]ds+C\Delta^3+\frac{3\Delta}{2P}\int_{t_{m}}^t\E[\Lambda_i\left(\bfX(s); k\right)]ds
+\frac{5}{2P}\int_{t_{m}}^t\E[\Lambda_i\left(\bfX(s);\sigma\right)]ds   .
\end{equation*}
Here, note that $\bfX(t)$ is independent of the random batches at $t_{m}$ so that
$$
\E\left[\E\left(\left|\chi_i\left(\bfX(s); k\right)\right|^2 |\cF_{m}\right)\right]=\E\left[\left|\chi_i\left(\bfX(s); k\right)\right|^2 \right]
\leq\frac{1}{P}\E[\Lambda_i\left(\bfX(s); k\right)].
$$
Finally, by iteration, one has
\begin{equation*}
 \E\left[\left|Z^i_t\right|^2\right]  \le C\int_{0}^t \E\left[\left|Z^i_{s}\right|^2\right]ds+C(m+1)\Delta^3+\frac{3\Delta}{2P}\int_{0}^t\E[\Lambda_i\left(\bfX(s); k\right)]ds
+\frac{5}{2P}\int_{0}^t\E[\Lambda_i\left(\bfX(s);\sigma\right)]ds   .
\end{equation*}

Applying Gr\"onwall's inequality (Theorem 1 in \cite{dragomir2003some}) yields the result as $(m+1)\Delta^3\le T\Delta^2$.
\end{proof}

\subsection{The nonlinear case}

In this subsection, we consider the analysis of the random batch approximation \eqref{eq:nonlinearrbm} for the nonlinear dependence case.

We will establish the approximation error estimate with the following assumption.
\begin{assumption}\label{ass:nonlinearlipschitz}
$f$ is one-sided Lipschitz. $A(\cdot):\R^r\to \R^d$, $k: \R^d\times\R^d\to \R^r$ and $\sigma: \R^d\to\R^{d\times m'}$ are all Lipschitz continuous and $\nabla^2 A$ is bounded.
\end{assumption}

Following similar argument as Proposition \ref{prpstn:boundedness2}, one can show that
\begin{prpstn}
Let Assumption \ref{ass:nonlinearlipschitz} hold and $\mu_0\in\cP_q$ with $q\ge 2$.
The process $\{X_t^i\}$ satisfies that
\begin{equation}
\sup _{0\leq t\leq T} \mathbb{E}[\left|X^{i}_{t}\right|^{q}]=\sup _{0\leq t\leq T} \mathbb{E}[\left|X^{1}_{t}\right|^{q}] \leq C_{q,T}, \quad \forall i.
\end{equation}
For the process $\tilde{X}_t^i$, conditioning on the sequence of random batches, one has
\begin{equation}
\sup_{0\le t\le T}\max_{1\le j\le N}\E[|\tilde{X}^j_t|^q  | \cF ]
\le C_{q,T},
\end{equation}
where $C_{q,T}$ is independent of $N$ and the sequence of batches.
\end{prpstn}
The proof is essentially the same and we omit.
Moreover, the symmetry Lemma \ref{lmm:symm} also holds for this case.

Below, let us focus on the error process $Z^i_t=\tilde{X}^i_t -X^i_t$. We find that
\[
\begin{split}
dZ^i_t & =(f(\tilde{X}^i_t)-f(X^i_t))\,dt
+\left(A\bigg(\frac{1}{P-1}\sum_{j\in \mathcal{C}_{q_{m}(i)},j\neq i} k(\tilde{X}^i_t, \tilde{X}^j_t)\bigg)
- A\bigg(\frac{1}{P-1}\sum_{j\in \mathcal{C}_{q_{m}(i)},j\neq i} k(X^i_t, X^j_t)\bigg) \right)\,dt \\
&\quad +\left(A\bigg(\frac{1}{P-1}\sum_{j\in \mathcal{C}_{q_{m}(i)},j\neq i} k(X^i_t, X^j_t)\bigg) - A\bigg(\frac{1}{N-1}\sum_{j\neq i} k(X^i_t, X^j_t)\bigg)\right)\,dt
+(\sigma(\tilde{X}^i_t)-\sigma(X^i_t))\,dW^i_t.
\end{split}
\]

Denote $F_i(\bfX(t)):= \frac{1}{N-1}\sum\limits_{j\neq i} k(X^i_t, X^j_t)$.
The main term that brings the difference is thus
\begin{equation*}
\begin{split}
A\Big(\chi_i(\bfX(t); k)+F_i(\bfX(t))\Big)-A(F_i(\bfX(t)))
&=\chi_i(\bfX(t);k)\cdot \nabla A(F_i(\bfX(t))) \\
&\quad+\frac{1}{2}\chi_i(\bfX(t); k)\otimes \chi_i(\bfX(t); k) : \nabla^2 A(F_i(\bfX(t))
+\xi \chi_i).
\end{split}
\end{equation*}
where $\xi$ is a constant in $(0,1)$.\\
As one may expect, the linear term is fine and can be estimated similarly as in
\cite{jin2020,jin2021convergence}. The main difference is the
quadratic term.
In fact, this term will be of the order of the variance so that
it contributes a term with order $1/P^2$ to the mean square error
(thus $O(1/P)$ to the error).

We first consider the fourth moments of the random deviation:
\begin{lmm}\label{lmm:4momentoffluc}
There is a universal constant $C$ such that for $N\ge 5$ that
\[
\E [|\chi_i(\b{x}; k)|^4]\le C\frac{Q(x; k)}{P^2},
\]
where
\[
\begin{aligned}
& Q(x; k)=[M_1(x; k)^4
+M_2(x; k)M_1(x;k)^2+M_2(x; k)^2+M_3(x; k)M_1(x; k)+M_4(x; k)],\\
& M_q(x; k)=\frac{1}{N}\sum_j |k(x_i, x_j)|^q, \quad q=1,2,3,4.
\end{aligned}
\]
\end{lmm}
The proof of this lemma is tedious, and we provide its sketch in
Appendix \ref{app:proofoflemma}.

\begin{thrm}\label{nonlinearcase}
Consider the nonlinear model  \eqref{eq:particlenonlinear}  and recall $Z^i_t=\tilde{X}^i_t -X^i_t$. Let Assumption \ref{ass:nonlinearlipschitz} hold.  Then,
\begin{equation}
\sup_{0\le t\le T}\max_i \E [|Z^i_t|^2]\le C\left(\Delta^2+\frac{\Delta}{P}\Lambda(k)\right)
+C(\|\nabla^2 A\|)\frac{Q(k)}{P^2},
\end{equation}
where $Q(k)=\sup\limits_{0\le t\le T}\E [Q(\bfX(t); k)]<+\infty$.
\end{thrm}

\begin{proof}
By It\^o's formula, one has for $t\in [t_m, t_{m+1})$ that
\allowdisplaybreaks[4]
\begin{align}
\frac{d}{dt}&\E\left[|Z^i_t|^2 | \cF_{m}\right]
\notag\\&\le   \Bigg\{2 \E\left[Z^i_t \cdot \left(f(\tilde{X}^i_t)-f(X^i_t)\right) \Big| \cF_{m} \right]+2\left\|\nabla A\right\|_{\infty} \E\bigg[\frac{|Z^i_t|}{P-1}\sum_{j\in \mathcal{C}_{q_{m}(i)},j\neq i}
\left|k\left(\tilde{X}^i_t, \tilde{X}^j_t\right) - k\left(X^i_t, X^j_t\right) \right| \Big| \cF_{m} \bigg]\Bigg\}
\notag\\&\quad+
\Bigg\{ 2\E\Big[ \chi_i\left(\bfX(t); k\right)\cdot \nabla A\left(F_i(\bfX(t))\right)\cdot Z^i_t  | \cF_{m} \Big]+2\E\left( \chi_i \otimes \chi_i : \nabla^2 A\left(F_i(\bfX(t))
+\xi \chi_i\right)\cdot Z^i_t | \cF_{m} \right)\Bigg\}
\notag\\
& \quad+ \left\{\E\bigg[\tr\left(\left(\sigma(\tilde{X}^i_t)-\sigma(X^i_t)\right) \left(\sigma(\tilde{X}^i_t)-\sigma(X^i_t)\right)^{\scriptscriptstyle T}\right) \Big| \cF_{m}   \bigg]\right\}
\notag\\&=:S_1+S_2+S_3,\notag
\end{align}
By Assumption \ref{ass:nonlinearlipschitz}, $S_1+S_3$ can be simply controlled by
$$
 \E[|Z^i_t|^2 | \cF_{m}]+\E\bigg[\frac{1}{P-1}\sum_{j\in \mathcal{C}_{q_{m}(i)}, j\neq i}
|Z^j_t|^2 | \cF_{m}\bigg]
$$
Consider the $S_2$ term. The second term  in $S_2$ will be bounded simply by
\[
\E[|Z^i_t|^2 | \cF_m]+\E[|\chi_i(\bfX(t))|^4 | \cF_m].
\]
For the first term in $S_2$, noticing that
\[
 \E\Big[ \chi_i(\bfX(t); k)\cdot \nabla A(F_i(\bfX(t)))\cdot Z^i_{t_m}   \Big]=0 ,
\]
one may again break
$Z^i_t=Z^i_{t_{m}}+(Z^i_t-Z^i_{t_{m}})$.
The treatment for $Z^i_t-Z^i_{t_{m}}$ is tedious but is similar to the proof of Theorem \ref{thm:batch1}. We omit the details.
Eventually, one has
\begin{align}
  \E\left[|Z^i_t|^2|\cF_{m}\right] &\le \E\left[|Z^i_{t_{m}}|^2|\cF_{m}\right]+C\int_{t_{m}}^t \E\left[|Z^i_{s}|^2|\cF_{m}\right]ds
+C\int_{t_{m}}^t\E\left[\frac{1}{P-1}\sum_{j\in \mathcal{C}_{q_{m}(i)},j\neq i} |Z^j_{s}|^2 |\cF_{m}\right]ds\notag\\
&\quad+\int_{t_{m}}^t\chi_i\left(\bfX(s); k\right)\cdot \nabla A\left(F_i\left(\bfX(s)\right)\right)\cdot Z^i_{t_{m}}ds +C\Delta^3\notag\\
&\quad+C\frac{\Delta}{2}\int_{t_{m}}^t\E\left[|\chi_i\left(\bfX(s); k\right)|^2 |\cF_{m}\right]ds
+C(1+\Delta)\int_{t_{m}}^t\E\left[\left|  \chi_{i}\left(\bfX(s); \sigma\right) \right|^4 \Big| \cF_{m} \right]ds.\notag
\end{align}

Taking expectation, using Lemma \ref{lmm:consistencyofrbm},
the symmetry (Lemma \ref{lmm:symm} also holds for this case and the proof is the same)
and Lemma \ref{lmm:4momentoffluc}, one has
$$
\E[|Z^i_t|^2]\le \E[|Z^i_{t_{m}}|^2]+C\int_{t_{m}}^t \E[|Z^i_s|^2]ds
+C\Delta^3+C\int_{t_{m}}^t\frac{\Lambda(k)\Delta}{P}ds+C\int_{t_{m}}^t\left(\left\|\nabla^2 A\right\|\right)\frac{Q(k)}{P^2}ds.
$$
Then, similar as in the proof of Theorem \ref{thm:batch1}, one can do iteration first and apply the Gr\"onwall's inequality to obtain the final estimate.
\end{proof}

Clearly, this error introduced by the nonlinearity is smaller compared to the one introduced by the random batch approximation in the noise term.

\subsection{Discussion on more general cases}

In this section, we discuss more general cases. In fact, the distribution dependence in the previous section does not include the cases when several statistical quantities are involved in the coefficients. For example, if the coefficient depends on the variance, then both the first and second moments will be involved.
These cases could be given by
\[
a(x, \mu)=f(x)+A\left(\int k_1(x, y)\mu(dy), \cdots, \int k_p(x, y)\mu(dy)\right)
\]
and
\[
b(x, \mu)=B\left(\int \sigma_1(x, y)\mu(dy), \cdots, \int \sigma_q(x, y)\mu(dy)\right).
\]
Here, $A, B, k_i, \sigma_j$ are all assumed to be Lipschitz and have bounded derivatives. $f$ is allowed to be one-sided Lipschitz.

The random approximation for such systems can be similarly performed. There is no big difference for the error estimate, while the details could be tedious. One can find that the strong error is again like
\[
\E [|Z^i_t|^2] \le \frac{C(\sigma)}{P}+\frac{C(k)}{P^2}+\frac{C\Delta}{P}.
\]
In other words, if there is law dependence in the coefficient of the noise, the mean square error  would be like $\cO(1/P)$. The mean error arising from the nonlinearity in the drift is like $\cO(1/P^2)$. These do not vanish as the step size tends to zero. However, the estimates are independent of $N$ so they could also save time if the required $N$ is large. In the special case when the diffusion coefficient $b$ is independent of $\mu$ and the drift linearly depends on $\mu$, the RBM approximation has an error that vanishes in the $\Delta\to 0$ limit, which is exactly the one studied in \cite{jin2020}.

\section{Numerical example}\label{part5}
%\tcr{The computational cost in table is for fixed terminal time so the error may not be the same for different items. Consider the case where there is dependence in the diffusion coefficient. Then, the error is like $\sqrt{1/P}\sim \epsilon$ for a given tolerance. Hence, in this case, $P=O(1/\epsilon^2)$ is fixed and the cost each iteration is $O(N/\epsilon^2)$. If we also require $P=1/\Delta^{\beta}$, then adjusting $\beta$ is equivalent to adjust $\Delta$ so that $\Delta=O(\epsilon^{2/\beta})$. Choosing larger $\beta$ makes the number of iteration smaller. If the discretization error of the SDE for each subinterval is $\Delta^{\alpha}$, then choosing $\beta=2\alpha$ would be the optimal one. In practice, the convergence rate seems better than $\sqrt{1/P}$, and this is not fully understood by us so it might be interesting to inverstigate this further.}

 In this section, we perform some numerical experiments to verify the theoretical analysis above.
In many applications of McKean-Vlasov SDEs, the dependence on the distribution in the coefficients are simply through the expected value of some kernel with respect to the law, which is a simple example of the cooperative interaction. This might due to the fact that many systems have a tendency to relax toward the center of gravity of the ensemble. See for example the model in \cite{dawson1983critical}.  Due to this reason, we will mainly consider examples of this form in this section.

%\cite{sharrock2021}
%\begin{equation*}
%X(t)=X(0)+\int_{0}^{t}\left(\lambda_1 X(s)\left(\lambda_2-\left|X(s)\right|\right)+\E X(s)\right)d s+\int_{0}^{t}\left( \xi\left|X(s)\right|^{3 / 2} +\E X(s)\right)d W(s)
%\end{equation*}
\begin{xmpl}
For linear cases, we consider the following one-dimensional SDE
\begin{equation}\label{example:linear}
dX(t)=\left(\lambda_1 X(t)\left(\lambda_2-\left|X(t)\right|\right)+\int_{\mathbb{R}}\left( X(t)-y\right)\mu(dy)\right)d t+\left( \lambda_3\left|X(t)\right|^{3 / 2} +\int_{\mathbb{R}}\left( X(t)-y\right)\mu(dy)\right)d W(t)
\end{equation}
 where  $\lambda_1$, $\lambda_2$ and $\lambda_3$ are positive constants.
\end{xmpl}

The corresponding system of interacting particles is given by
\begin{equation}\label{example:particlesystem}
dX^{i,N}_{t}=\bigg(\lambda_1 X^{i,N}_{t}\left(\lambda_2-|X^{i,N}_{t}|\right)+\frac{1}{N-1}\sum_{j\neq i}^{N}(X^{i,N}_{t}-X^{j,N}_{t})\bigg)d t+\bigg( \lambda_3|X^{i,N}_{t}|^{3 / 2} +\frac{1}{N-1}\sum_{j\neq i}^{N}(X^{i,N}_{t}-X^{j,N}_{t})\bigg)d W_t^i.
\end{equation}

\begin{figure}[htbp]
\centering
\includegraphics[width=3in]{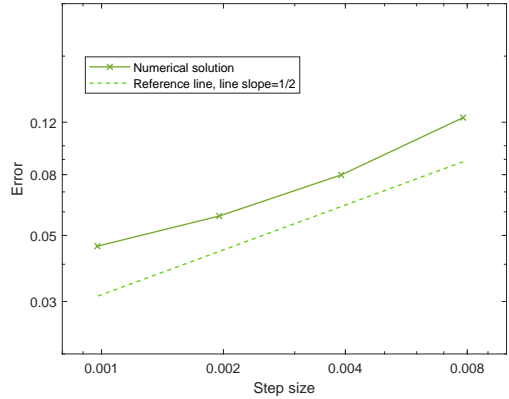}
\caption{The convergence order of the truncated EM method applied to the interacting particle system related to (\ref{example:linear}) without using the random batch technique. \label{figure:linear1}}
\end{figure}

%In what follows, we set parameters $X_{0}=1, \lambda_1=2.5, \lambda_2=1, \lambda_3=1$. We choose particle number $N=2^{12}$ and obtain \tcb{1000} sample paths (i.e., repeat the simulation for $1000$ times) to calculate the mean square error at $T=1$ with four different  step sizes: $\Delta=2^{p}$ for $-10 \leq p \leq -7$.
%\tcr{This is an absolute error $\sqrt{\mathbb{E}[|X_{t}^{i, N}-X_{t}^{i, \Delta}|^2]}$, where the expectation is approximated by the Monte Carlo using all the particles for all the sample paths.}
In the subsequent analysis, we fix the following parameters: $X_{0}^{i}=1, \lambda_1=2.5, \lambda_2=1, \lambda_3=1$, for all $i\in\{1,\ldots,N\}$. We select a particle number of $N=2^{12}$ and generate a total of 1000 sample paths (i.e., conducting the simulation $1000$ times) to evaluate the mean square error at $T=1$. This evaluation is performed using four different step sizes: $\Delta=2^{p}$ for $-10 \leq p \leq -7$.
The mean square error, denoted by $\sqrt{\mathbb{E}[|X_{t}^{i, N}-X_{t}^{i, \Delta}|^2]}$, is calculated by approximating the expectation via Monte Carlo simulations, utilizing all the particles for all the sample paths.

Figure \ref{figure:linear1} displays the resulting log-log error plot along with a reference line for slope one-half. The slopes of the two curves seem to line up nicely, which confirms the theoretical result in Theorem \ref{thm3.4} when the truncated EM method applied to the interacting particle system corresponding to (\ref{example:linear}) directly. Here we choose the tame Euler scheme with time step size $\Delta=2^{-12}$ and the same number of particles $N=2^{12}$ for producing a reference solution.

To validate Theorem \ref{thm:batch1}, we also  employ  truncated type methods together with the random batch technique to do a simulation.
From Corollary \ref{crll}, we note that the order of convergence of our approach to approximating the interacting particle system related to the McKean-Vlasov SDE (\ref{example:linear}) is at least $\beta/2$.

\begin{figure}[htbp]
\centering
\includegraphics[width=3in]{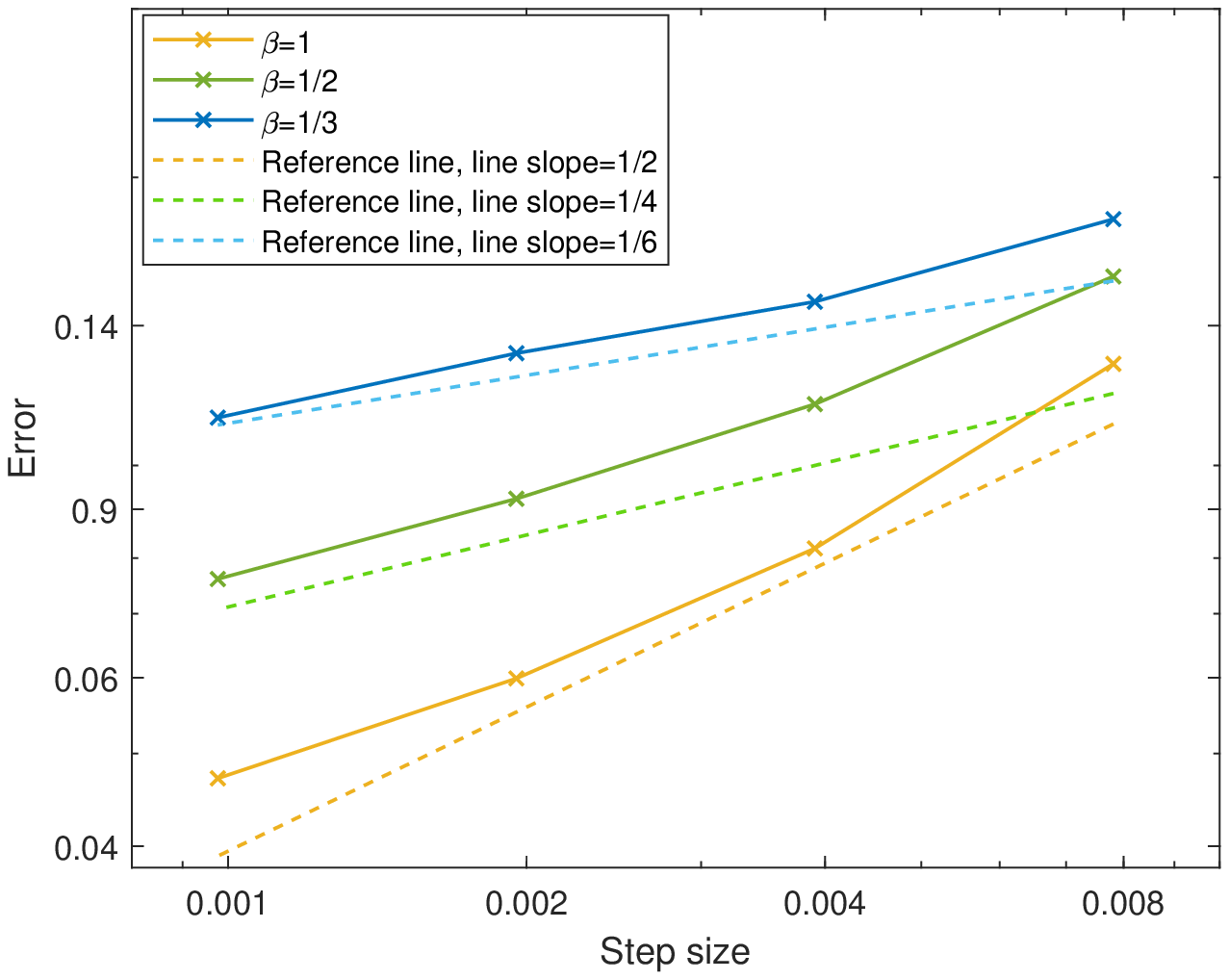}
    \caption{Convergence order of truncated EM method applied to the interacting particle system related to (\ref{example:linear}) } with random batch technique for $\beta=1,1/2,1/3$.\label{figure:linear3}
\end{figure}

Now we generate the reference solution using the truncated Euler method \eqref{INn} with a step size of $\Delta=2^{-12}$. Subsequently, with three different values of $\beta$ chosen, we obtain three distinct batch sizes denoted by $P=1/\Delta^\beta$. The numerical errors associated with these batch sizes are plotted in Figure \ref{figure:linear3}.
%To determine the convergence order, we employ the least squares fit method \cite{higham2001}, resulting in convergence order values of \tcr{$0.4766$, $0.3477$, and $0.2247$} for $\beta=1$, $1/2$, and $1/3$, respectively.
Hence, our findings align with the expected strong order of convergence stated in Corollary \ref{crll}.%\tcr{Theorem \ref{thm:batch1}(Corollary \ref{crll})}.

%\begin{figure}[htbp]
%\centering
%\includegraphics[width=3in]{convergence}
%    \caption{Convergence order of truncated EM method with random batch technique for $\beta=1,1/2,1/3$.\label{vs}}
%\end{figure}

 Consider that the main advantage of random batch methods is to save the computational cost, we compare the  computation time of the truncated Euler method using a random batch technique and the counterpart without using the random batch technique. In Table \ref{table:1}, we denote the former by the notation `TEMwRBM', and the notation `TEM' represents the latter. Let us consider SDE (\ref{example:linear}) where there is a dependence in the diffusion coefficient. It is easy to find that the execution time significantly
enlarges when increasing the number of particles. Here we fix the time step size $\Delta=2^{-7}$. The `TEM' and the `TEMwRBM' run 1000 sample paths to obtain the average time consumption respectively. %The results in Table \ref{table:1} demonstrate better efficiency with a larger amount of the particle number.

Upon implementing the random batch technique, we can approximate the error, as stated in (\ref{eq:mainestimate}), using $\sqrt{1/P} \sim \epsilon$ for a given tolerance $\epsilon$. In this scenario, $P=O(1/\epsilon^2)$ remains fixed, resulting in a cost per iteration of $O(NP)$. However, when $P=1/\Delta^{\beta}$ is also required, the cost per iteration increases as the value of $\beta$ grows.
		On the other hand, the computational cost presented in Table \ref{table:1} corresponds to a fixed terminal time. Notably, adjusting $\beta$ is equivalent to modifying $\Delta$ such that $\Delta=O(\epsilon^{2/\beta})$. Increasing the value of $\beta$ results in a larger step size, consequently reducing the number of iterations. Thus, an optimal value of $\beta$ exists in terms of the overall computational cost. However, we have yet to fully comprehend this phenomenon, and further investigation would be both intriguing and worthwhile.%When the discretization error of the SDE for each subinterval is $\Delta^{\alpha}$, selecting $\beta=2\alpha$ is considered optimal. Interestingly, in practice, the convergence rate appears to be better than $\sqrt{1/P}$. However, we have not fully understood this phenomenon, and further investigation would be intriguing.

\begin{table*}[t]
\begin{center}
\caption{Execution time comparison}
\label{table:1}
\setlength{\tabcolsep}{7mm}{
\begin{tabular}{|c|r|r|r|r|}
\hline
\multirow{2}{*}{N}&\multirow{2}{*}{TEM}&\multicolumn{3}{c|}{TEMwRBM}\\
\cline{3-5}
 & &  $P=1/\Delta$ &  $P=  1/\sqrt{\Delta}$&$P=  1/\sqrt[3]{\Delta}$\\
\hline
$2^{12}$& 4.8672s&1.3771s&0.9657s&0.9559s\\
\hline
$2^{14}$& 84.1268s&5.6311s&3.9629s&4.1012s\\
\hline
$2^{16}$ & 1143.8016s&21.7272s&14.9615s&15.4843s\\
\hline
\end{tabular}}
\end{center}
\end{table*}

\begin{xmpl}
When there is no interaction in the diffusion part,  (\ref{example:linear}) degenerates to
\begin{equation}\label{example:linear2}
dX(t)=\left(\lambda_1 X(t)\left(\lambda_2-\left|X(t)\right|\right)+\int_{\mathbb{R}}\left( X(t)-y\right)\mu(dy)\right)d t+ \lambda_3\left|X(t)\right|^{3 / 2} d W(t).
\end{equation}
\end{xmpl}

The corresponding system of interacting particles is given by
\begin{equation}\label{example:particlesystem2}
dX^{i,N}_{t}=\bigg(\lambda_1 X^{i,N}_{t}\left(\lambda_2-|X^{i,N}_{t}|\right)+\frac{1}{N-1}\sum_{j\neq i}^{N}(X^{i,N}_{t}-X^{j,N}_{t})\bigg)d t+ \lambda_3|X^{i,N}_{t}|^{3 / 2} d W_t^i.
\end{equation}

In what follows, we set parameters $X_{0}^{i}=1, \lambda_1=2.5, \lambda_2=1, \lambda_3=1$, for all $i\in\{1,\ldots,N\}$. We choose particle number $N=2^{12}$ and obtain 200 sample paths (i.e., repeat the simulation for $200$ times) to calculate the mean square error at $T=1$ with four different  step sizes: $\Delta=2^{p}$ for $-10 \leq p \leq -8$.

\begin{figure}[htbp]
\centering
\includegraphics[width=3in]{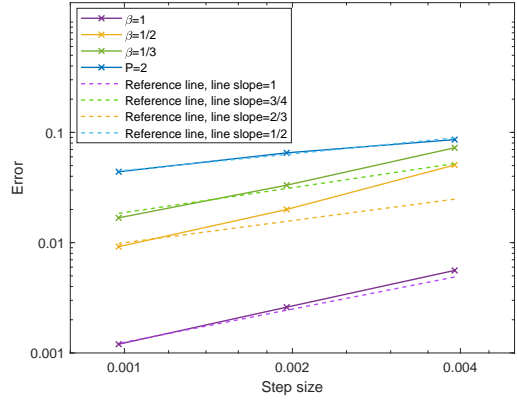}
\caption{The convergence order of truncated Milstein method with a random batch size  $P=1/\Delta^{\beta}$  with $\beta=1, 1/2, 1/3$ and $P=2$ applied to the interactive particle system corresponding to \eqref{example:linear2}. \label{figure:linear2} }
\end{figure}

The conclusion (\ref{eq:mainestimate}) implies that we can obtain the first-order convergence when approximating the interacting particle system associated to (\ref{example:linear2}) by a random batch technique with $\beta=1$ (the batch size $P=1/\Delta$). At the same time strong convergence order of a numerical method  should be one. Therefore we employ the truncated Milstein scheme \cite{guo2018} to verify  the first result of the Corollary \ref{crll}. The numerical results illustrated in Figure \ref{figure:linear2} confirm our theoretical result, where a dashed red line with a slope of one is provided for reference.  Moreover, we conducted a comparison of the results with varying batch sizes. When considering $P=1/\Delta^{\beta}$ with $\beta=1, 1/2, 1/3$, as well as $P=2$, as shown in Figure \ref{figure:linear2} and agrees with the theory in Corollary \ref{crll}. 

\begin{xmpl}
For nonlinear cases, we consider the following SDE
\begin{equation*}\label{example:nonlinear}
dX(t)=\left(\lambda_1 X(t)\left(\lambda_2-\left|X(t)\right|\right)+\sin\Big(\int_{\mathbb{R}}\left( X(t)-y\right)\mu(dy)\Big)\right)d t
 +\lambda_3\left|X(t)\right|^{3 / 2} d W(t)
\end{equation*}
where  $\lambda_1$, $\lambda_2$ and $\lambda_3$ are positive constants.
\end{xmpl}

The corresponding system of interacting particles is given by
\begin{equation*}\label{example:nonlinearparticlesystem}
dX^{i,N}_{t}=\bigg(\lambda_1 X^{i,N}_{t}\left(\lambda_2-|X^{i,N}_{t}|\right)+\sin\Big(\frac{1}{N-1}\sum_{j\neq i}^{N}(X^{i,N}_{t}-X^{j,N}_{t})\Big)\bigg)d t
+\lambda_3|X^{i,N}_{t}|^{3 / 2} d W_t^i.
\end{equation*}

\begin{figure}[htbp]
\centering
\includegraphics[width=3in]{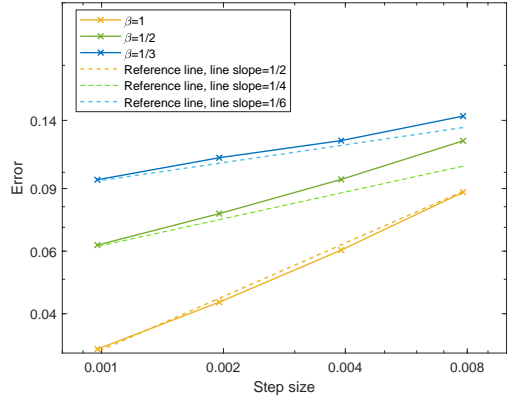}
    \caption{The convergence order of the truncated EM method and random batch technique with $\beta=1,1/2,1/3 $ in nonlinear cases.\label{figure:nonlinear}}
\end{figure}

We examine three distinct batch sizes represented by $P=1/\Delta^\beta$. The numerical errors are illustrated in Figure \ref{figure:nonlinear}, where three dashed red lines with slopes of 1/2, 1/4, and 1/6 are provided as reference. %By performing a least squares fit of the errors, we obtain convergence order values of \tcr{$0.4887$, $0.3242$, and $0.1946$} for $\beta=1$, $1/2$, and $1/3$, respectively.
These results align with the expected strong order of convergence stated in Theorem \ref{nonlinearcase}.
%\tcr{We have three kinds of batch size $P=1/\Delta^\beta$, then the numerical errors are plotted in Figure \ref{figure:nonlinear},  where a dashed red line with a slope of one-half is provided for reference. By least squares fit of error we have the value of convergence order are \tcr{$0.4887$, $0.3242$ and $0.1946$} for $\beta=1$, $1/2$, and $1/3$ respectively.  Hence, our results are consistent with a strong order of convergence from Theorem \ref{nonlinearcase}.}

%\begin{figure}[!b]
%\centering
%\includegraphics{02}
%\caption{This is an example for appendix figure\label{fig4}}
%\end{figure}

%\begin{table}[t]
%%\begin{center}
%\caption{Execution time comparison\label{biao1}}
%\begin{tabular*}{\textwidth}{@{\extracolsep\fill}lcccc@{}}
%\toprule
%\multirow{2}{*}{U}&\multirow{2}{*}{TEM}&\multicolumn{3}{c}{TEMwRBM}\\
%\cline{3-5}
% & &  $P=1/\Delta$ &  $P=  1/\sqrt{\Delta}$&$P=  1/\sqrt[3]{\Delta}$\\
%\midrule
%$2^{12}$& 15.89s&11.80s&3.05s&3.02s\\
%$2^{14}$& 239.99s&56.03s&12.06s&11.66s\\
%$2^{16}$ & 2529.15s&211.71s&49.18s&50.26s\\
%\botrule
%\end{tabular*}
%%\end{center}
%\end{table}

%%%%%%%%%%%%%%

%The conclusion (\ref{eq:mainestimate}) implies that we can obtain the first-order convergence when approximating the interacting particle system associated to (\ref{example:linear2}) by a random batch technique with $\beta=1$ (the batch size $P=1/\Delta$). At the same time strong convergence order of a numerical method  should be one. Therefore we employ the truncated Milstein scheme \cite{guo2018} to verify  \tcb{Corollary \ref{crll}}. The numerical results illustrated in Figure \ref{figure:linear2} confirm our theoretical result, where a dashed red line with a slope of one is provided for reference.

\appendix

\section{Proof of Proposition \ref{pro:momentbound1}}\label{app:bound1}

\begin{proof}
Here, we only illustrate the moments control for $X_t^i$. The one for $\bar{X}_t^i$ is similar and simpler.

For \eqref{I},  using It\^o’s formula,
\begin{equation*}
%\begin{aligned}
\mathbb{E}[|X_{t}^{i}|^{q}]\leq  \mathbb{E}[|X_{0}^{i}|^{q}]+q \mathbb{E}\int_{0}^{t}|X_{s}^{i}|^{q-2}\langle X_{s}^{i}, a\left(X_{s}^{i}, \mu_{s}^{X}\right)\rangle ds+\frac{q(q-1)}{2} \mathbb{E}\int_{0}^{t}|X_{s}^{i}|^{q-2}\|b\left(X_{s}^{i}, \mu_{s}^{X}\right)\|^2 ds.
%\end{aligned}
\end{equation*}
By \eqref{Khasminskii1}, we get
\begin{equation*}
\mathbb{E}[|X_{t}^{i}|^{q}]
 \leq  \mathbb{E}[|X_{0}^{i}|^{q}]+  qL \mathbb{E}\int_{0}^{t} |X_{s}^{i}|^{q-2}[(1+|X_{s}^{i}|)^{2}+\mathcal{W}_{2}^{2}\left(\mu_s^X, \delta_0\right)]ds.
\end{equation*}
Apply the Young's inequality and Lemma \ref{lmm:w2aux},
\begin{equation*}
\mathbb{E}[|X_{t}^{i}|^{q}] \leq  C(1+\mathbb{E}[|X_{0}^{i}|^{q}])+ C  \int_{0}^{t}\mathbb{E}[|X_{s}^{i}|^{q}]ds+C \int_{0}^{t}\mathbb{E}\left(\frac{1}{N}\sum_{\tcb{j=1}}^N|X_{s}^{j}|^2\right)^{q/2}ds.
\end{equation*}
For the last term, we apply  the Minkowski inequality and the symmetry consecutively,
\[
\mathbb{E}\left(\frac{1}{N} \sum_{j=1}^{N}\left|X_{s}^{j}\right|^{2}\right)^{q/2}
\le \left(\frac{1}{N}\sum_{j=1}^N \|\left|X_{s}^{j}\right|^{2}\|_{L^{q/2}} \right)^{q/2}
=\E[\left|X_{s}^{i}\right|^{q}].
\]

Hence,
\begin{equation*}
\begin{aligned}
\mathbb{E}[|X_{t}^{i}|^{q}]& \leq  C  \int_{0}^{t}\mathbb{E}[|X_{s}^{i}|^{q}]ds+C(1+\mathbb{E}[|X_{0}^{i}|^{q}]),
\end{aligned}
\end{equation*}
and applying the Gr\"onwall’s inequality finishes the estimate.
\end{proof}

\section{Proof of the auxilliary lemma \ref{lmm:4momentoffluc}}\label{app:proofoflemma}

Recall that $q(i)$ means the index of the batch that particle $i$ falls into.
Introduce the indicator function
\[
I_{ij}=\begin{cases}
1 &  q(i)=q(j),\\
0 & \text{~otherwise}.
\end{cases}
\]
Then, we have the simple lemma:
\begin{lmm}\label{lmm:indicatorexpec}
Let $q<N$. Let $i, j_{j'}$, $1\le j' \le q$ be all distinct. Then,
\[
\E \prod_{j'=1}^q I_{i j_{j'}}=\prod_{j'=1}^q \frac{P-j'}{N-j'}.
\]
\end{lmm}
\begin{proof}
If $P \le q$, clearly both sides are zero so the inequality holds.
Now, assume $P>q$.

The event for which $\prod\limits_{j'=1}^q I_{i j_{j'}}=1$ is when all the particles are in the same batch. The proof would be similar to the argument in \cite[Lemma 3.1]{jin2020}.
There are
\[
M(n):=\frac{(nP)!}{(P!)^n n!}
\]
ways to divide $N=nP$ particles into $n$ batches, without consideration of the order of batches. Then, the probability that $i, j_1,\cdots, j_q$ are in the same batch is given by
\[
\frac{{{nP-q-1}\choose{P-q-1}}M(n-1)}{M(n)}=\prod_{j'=1}^q \frac{P-j'}{N-j'}.
\]
\end{proof}

\begin{proof}[Proof of Lemma \ref{lmm:4momentoffluc}]
Clearly,
\[
\chi_i(\b{x}; k)=\sum_{j=1, j\neq i}^N \left(\frac{1}{P-1}I_{ij}-\frac{1}{N-1}\right)k(x_i, x_j).
\]
It follows that
\[
|\chi_i(\b{x}; k)|^4
=\sum_{j_1 j_2 j_3 j_4,\neq i}\zeta(j_1, j_2, j_3, j_4) k(x_i, x_{j_1})\cdot k(x_i, x_{j_2})
k(x_i, x_{j_3})\cdot k(x_i, x_{j_4}).
\]
Here, for matrices, $k(x_i, x_{j_1})\cdot k(x_i, x_{j_2})$ means $\tr(k(x_i, x_{j_1}) k(x_i, x_{j_2})^T)$. The random variable $\zeta$ is given by
\begin{equation*}
\begin{split}
\zeta(j_1, j_2, j_3, j_4)
=&\frac{1}{(P-1)^4}\prod_{j'}I_{i,j_{j'}}-\frac{1}{(P-1)^3(N-1)}
(I_{ij_1}I_{ij_2}I_{ij_3}+I_{ij_1}I_{ij_2}I_{ij_4}
+I_{ij_1}I_{ij_3}I_{ij_4}+I_{ij_2}I_{ij_3}I_{ij_4})\\
+&\frac{1}{(P-1)^2(N-1)^2}(I_{ij_1}I_{ij_2}+I_{ij_1}I_{ij_3}
+I_{ij_1}I_{ij_4}+I_{ij_2}I_{ij_3}+I_{ij_2}I_{ij_4}+I_{ij_3}I_{ij_4})\\
-&\frac{1}{(P-1)(N-1)^3}(I_{i j_1}+I_{i j_2}+I_{i j_3}+I_{i j_4})
+\frac{1}{(N-1)^4}.
\end{split}
\end{equation*}
Here, we consider several cases.

{\bf Case 1}: $j_1, j_2,j_3, j_4$ are all distinct. There are
$(N-1)(N-2)(N-3)(N-4)$ such terms. By Lemma \ref{lmm:indicatorexpec}, one has
\begin{equation*}
\begin{split}
\E \zeta(j_1, j_2, j_3, j_4)
=&~\frac{(P-2)(P-3)(P-4)}{(P-1)^3 \prod_{j'=1}^4 (N-j')}
-\frac{4(P-2)(P-3)}{(P-1)^2(N-1)^2(N-2)(N-3)}
\\&+\frac{6(P-2)}{(P-1)(N-1)^3(N-2)}-\frac{3}{(N-1)^4}.
\end{split}
\end{equation*}
Adding these up, one has $Num/[(P-1)^3(N-1)^4(N-2)(N-3)(N-4)]$. We need to check the coefficients of $P^3 N^3, P^2 N^3$ and $P^3N^2$
in the numerator $Num$. After some calculation, one finds that these are all zero. Hence,
\[
\E \zeta(j_1, j_2, j_3, j_4) \le C\frac{PN^3+P^2N^2+P^3N}{P^3 N^7}
\le C\frac{1}{P^2N^4}.
\]
It is then easy to see that these terms are controlled by
$C\frac{1}{P^2}M_1(x; k)^4$.

{\bf Case 2}: There are two indices that are the same. There are $6(N-1)(N-2)(N-3)$ such terms.
\begin{multline*}
\E \zeta(j_1, j_2, j_3, j_4)
=\frac{1}{(P-1)^4}\frac{(P-1)(P-2)(P-3)}{(N-1)(N-2)(N-3)}\\
-\frac{1}{(P-1)^3(N-1)}\left(2\frac{(P-1)(P-2)}{(N-1)(N-2)}
+2\frac{(P-1)(P-2)(P-3)}{(N-1)(N-2)(N-3)}\right)\\
+\frac{1}{(P-1)^2(N-1)^2}\left(\frac{P-1}{N-1}+5\frac{(P-1)(P-2)}{(N-1)(N-2)}\right)
-\frac{3}{(N-1)^4}.
\end{multline*}
Similarly, by combining terms and checking the corresponding coefficients
for some big terms, one can find that this is controlled by
\[
C\frac{PN^2+P^2N}{P^3 N^5}+C\frac{P+1}{P^2N^4}
\le C\frac{1}{P^2N^3}.
\]
Hence, this case is controlled by $C\frac{1}{P^2}M_1^2 M_2$.

{\bf Case 3} Two indices are the same while the other two are the same.
There are $3 (N-1)(N-2)$ such terms.
\begin{multline*}
\E \zeta(j_1, j_2, j_3, j_4)
=\frac{1}{(P-1)^4}\frac{(P-1)(P-2)}{(N-1)(N-2)}
-\frac{4}{(P-1)^3(N-1)}\frac{(P-1)(P-2)}{(N-1)(N-2)}\\
+O(\frac{1}{(P-1)^2(N-1)^2}+\frac{1}{(P-1)(N-1)^3}).
\end{multline*}
This is clearly controlled by $C/P^2$, and thus the sum of such terms
is controlled by $M_2^2/P^2$.

{\bf Case 4} Three indices are the same. There are $4(N-1)(N-2)$ such terms.

This situation can be treated similarly as Case 3. This is controlled by $M_3M_1/P^2$.

{\bf Case 5} All indices are the same. There are $N-1$ such terms.
This is controlled by $M_4/P^3$.

These then add up to the estimate in the lemma.
\end{proof}

\section*{Acknowledgments}
%The authors thank the anonymous reviewers for their valuable suggestions. This work is supported in part by funds from the National Science Foundation (NSF: \# 1636933 and \# 1920920).

This work was partially supported by the National Key R\&D Program of China, Project Number 2021YFA1002800.
The work of Q. Guo and J. He were partially supported by the National Natural Science Foundation of China (12271368, 11871343)  and Shanghai Rising-Star Program (22QA1406900).
The work of L. Li was partially sponsored by the Strategic Priority Research Program of Chinese Academy of Sciences, Grant No. XDA25010403 and NSFC Grant No. 11901389, 12031013, and Shanghai Science and Technology Commission Grant No. 20JC144100, 21JC1403700.
The authors would like to thank Professor Shi Jin for his valuable comments.

%We would like to thank David \v{S}i\v{s}ka and Lukasz Szpruch, both from University of Edinburgh, for allowing us to use the above mentioned Lemma~\ref{lem:DL} from their working paper which played a crucial role in the proof of Theorem \ref{thm:eu}.

All authors contributed equally to the manuscript.

\bibliographystyle{plain}
\bibliography{sdereference}

%%-----------------------------
%%      your bibliography
%%-----------------------------
\end{document}